\documentclass[11pt]{article}
\usepackage{setspace}
\usepackage{amsmath} 
\usepackage{amssymb}
\usepackage{latexsym}
\usepackage{xcolor}
\usepackage{fancyhdr}
\allowdisplaybreaks 

 \numberwithin{equation}{section}
\newtheorem{thm}[equation]{Theorem}
\newtheorem{prop}[equation]{Proposition}
\newtheorem{defn}[equation]{Definition}
\newtheorem{rem}[equation]{Remark}
\newtheorem{lem}[equation]{Lemma}

\usepackage{xcolor}

\title{Classes of kernels and continuity properties of the double layer potential in H\"{o}lder spaces}

\author{  
Massimo Lanza de Cristoforis
\\
Dipartimento di Matematica `Tullio Levi-Civita', 
\\
Universit\`a degli Studi di Padova, 
\\
Via Trieste 63, Padova 35121, 
Italy. 
\\
E-mail: mldc@math.unipd.it   }

 \date{\ }

\begin{document} 

\maketitle

\noindent
{\bf Abstract:}  We prove the validity of regularizing properties of the boundary integral operator corresponding to the double layer potential associated to the fundamental solution of a {\em nonhomogeneous} second order elliptic differential operator with constant coefficients in H\"{o}lder spaces by exploiting an estimate on the maximal function  of the tangential gradient with respect to the first variable of the kernel of the double layer potential and  by exploiting specific imbedding and multiplication properties in certain classes of integral operators and  a generalization of a result for integral operators on differentiable manifolds. 
  \vspace{\baselineskip}

\noindent
{\bf Keywords:} 
Double layer potential, second order differential operators with constant
coefficients, boundary behavior, H\"{o}lder spaces.\par

\noindent   
{{\bf 2020 Mathematics Subject Classification:}}   31B10

 \section{Introduction} In this paper, we consider the double layer potential associated to the fundamental solution of a second order differential operator with constant coefficients in H\"{o}lder spaces. Unless otherwise specified,  we assume  throughout the paper that
\[
n\in {\mathbb{N}}\setminus\{0,1\}\,,
\]
where ${\mathbb{N}}$ denotes the set of natural numbers including $0$. Let $\alpha\in[0,1]$, $m\in {\mathbb{N}}\setminus\{0\}$. Let $\Omega$ be a bounded open subset of ${\mathbb{R}}^{n}$ of class $C^{m,\alpha}$ and we understand that $C^{m,0}\equiv C^{m}$. For the definition and properties of the classical  Schauder spaces we refer for example to \cite[Chap.~2]{DaLaMu21},  \cite[\S 2]{DoLa17}. We employ the same notation of reference \cite{DoLa17} with Dondi that we now introduce. 

Let $\nu_\Omega$ or simply 
  $\nu \equiv (\nu_{l})_{l=1,\dots,n}$ denote the external unit normal to $\partial\Omega$.  Let $N_{2}$ denote the number of multi-indexes $\gamma\in {\mathbb{N}}^{n}$ with $|\gamma|\leq 2$. For each 
\begin{equation}
\label{introd0}
{\mathbf{a}}\equiv (a_{\gamma})_{|\gamma|\leq 2}\in {\mathbb{C}}^{N_{2}}\,, 
\end{equation}
we set 
\[
a^{(2)}\equiv (a_{lj} )_{l,j=1,\dots,n}\qquad
a^{(1)}\equiv (a_{j})_{j=1,\dots,n}\qquad
a\equiv a_{0}\,.
\]
with $a_{lj} \equiv 2^{-1}a_{e_{l}+e_{j}}$ for $j\neq l$, $a_{jj} \equiv
 a_{e_{j}+e_{j}}$,
and $a_{j}\equiv a_{e_{j}}$, where $\{e_{j}:\,j=1,\dots,n\}$  is the canonical basis of ${\mathbb{R}}^{n}$. We note that the matrix $a^{(2)}$ is symmetric. 
Then we assume that 
  ${\mathbf{a}}\in  {\mathbb{C}}^{N_{2}}$ satisfies the following ellipticity assumption
\begin{equation}
\label{ellip}
\inf_{
\xi\in {\mathbb{R}}^{n}, |\xi|=1
}{\mathrm{Re}}\,\left\{
 \sum_{|\gamma|=2}a_{\gamma}\xi^{\gamma}\right\} >0\,,
\end{equation}
and we consider  the case in which
\begin{equation}
\label{symr}
a_{lj} \in {\mathbb{R}}\qquad\forall  l,j=1,\dots,n\,.
\end{equation}
Then we introduce the operators
\begin{eqnarray*}
P[{\mathbf{a}},D]u&\equiv&\sum_{l,j=1}^{n}\partial_{x_{l}}(a_{lj}\partial_{x_{j}}u)
+
\sum_{l=1}^{n}a_{l}\partial_{x_{l}}u+au\,,
\\
B_{\Omega}^{*}v&\equiv&\sum_{l,j=1}^{n} \overline{a}_{jl}\nu_{l}\partial_{x_{j}}v
-\sum_{l=1}^{n}\nu_{l}\overline{a}_{l}v\,,
\end{eqnarray*}
for all $u,v\in C^{2}(\overline{\Omega})$, and a fundamental solution $S_{{\mathbf{a}} }$ of $P[{\mathbf{a}},D]$, and the  boundary integral operator corresponding to the  double layer potential 
\begin{eqnarray}
\label{introd3}
\lefteqn{
W_\Omega[{\mathbf{a}},S_{{\mathbf{a}}}   ,\mu](x) \equiv 
\int_{\partial\Omega}\mu (y)\overline{B^{*}_{\Omega,y}}\left(S_{{\mathbf{a}}}(x-y)\right)
\,d\sigma_{y}
}
\\  \nonumber
&&
\qquad
=-\int_{\partial\Omega}\mu(y)\sum_{l,j=1}^{n} a_{jl}\nu_{l}(y)\frac{\partial S_{ {\mathbf{a}} } }{\partial x_{j}}(x-y)\,d\sigma_{y}
\\  \nonumber
&&
\qquad\quad
-\int_{\partial\Omega}\mu(y)\sum_{l=1}^{n}\nu_{l}(y)a_{l}
S_{ {\mathbf{a}} }(x-y)\,d\sigma_{y} \qquad\forall x\in \partial\Omega\,,
\end{eqnarray}
where the density or moment $\mu$ is a function  from $\partial\Omega$ to ${\mathbb{C}}$. Here the subscript $y$ of $\overline{B^{*}_{\Omega,y}}$ means that we are taking $y$ as variable of the differential operator $\overline{B^{*}_{\Omega,y}}$. 
The role of the double layer potential in the solution of boundary value problems for the operator $P[{\mathbf{a}},D]$ is well known (cf.~\textit{e.g.}, 
G\"{u}nter~\cite{Gu67}, Kupradze,  Gegelia,  Basheleishvili and 
 Burchuladze~\cite{KuGeBaBu79}, Mikhlin \cite{Mik70}.) 
 
 The analysis of the continuity and compactness properties of $W_\Omega[{\mathbf{a}},S_{{\mathbf{a}}},\cdot]$ is  a classical topic and several results in the literature show that  $W_\Omega[{\mathbf{a}},S_{{\mathbf{a}}},\cdot]$ improves the regularity of  H\"{o}lder continuous functions on $\partial\Omega$. We briefly recall some references (see also \cite{DoLa17}).
 
 In case $n=3$ and $\Omega$ is of class $C^{1,\alpha}$ and 
$S_{{\mathbf{a}}}$ is the fundamental solution of the Laplace operator,
it has long been known that  $W_\Omega[{\mathbf{a}},S_{{\mathbf{a}}},\cdot]$
is  a linear and compact operator in $C^{1,\alpha}(\partial\Omega)$ and is linear and continuous from $C^{0 }(\partial\Omega)$ to $C^{0,\alpha}(\partial\Omega)$ (cf.~
Schauder~\cite{Sc31}, \cite{Sc32}, Miranda~\cite{Mi65}.)

In case  $n=3$, $m\geq 2$ and $\Omega$ is of class $C^{m,\alpha}$ and if  $P[{\mathbf{a}},D]$ is the  Laplace operator,  G\"{u}nter~\cite[Appendix, \S\ IV, Thm.~3]{Gu67} has 
proved that   $W[\partial\Omega ,{\mathbf{a}},S_{{\mathbf{a}}},\cdot]$ is bounded from $C^{m-2,\alpha}(\partial\Omega)$ to $C^{m-1,\alpha'}(\partial\Omega)$ for $\alpha'\in]0,\alpha[$.\par
	
In case $n\geq 2$, $\alpha\in]0,1]$,  O.~Chkadua \cite{Chk23} has pointed out  that one could exploit Kupradze,  Gegelia,  Basheleishvili and 
 Burchuladze~\cite[Chap.~IV, Sect.~2, Thm 2.9, Chap. IV, Sect.~3, Theorems 3.26 and  3.28]{KuGeBaBu79} and  prove that if $\Omega$ is of class $C^{m,\alpha}$, then $W[\partial\Omega ,{\mathbf{a}},S_{{\mathbf{a}}},\cdot]$ is bounded from $C^{m-1,\alpha'}(\partial\Omega)$ to $C^{m,\alpha'}(\partial\Omega)$ for $\alpha'\in]0,\alpha[$.\par

In case $n=3$  and $\Omega$ is of class $C^{2}$   and if  $P[{\mathbf{a}},D]$ is the  Helmholtz operator, Colton and Kress~\cite{CoKr83} have 
developed previous work of G\"{u}nter~\cite{Gu67} and Mikhlin~\cite{Mik70}  and proved that the operator $W_\Omega[{\mathbf{a}},S_{{\mathbf{a}}},\cdot]$
is bounded from $C^{0,\alpha}(\partial\Omega)$ to $C^{1,\alpha}(\partial\Omega)$ and that accordingly it
 is  compact in $C^{1,\alpha}(\partial\Omega)$.
 
  In case $n\geq 2$, $\alpha\in]0,1[$  and $\Omega$ is of class $C^{2}$   and if  $P[{\mathbf{a}},D]$ is the  Laplace operator, 
Hsiao and Wendland \cite[Remark 1.2.1]{HsWe08} 
deduce  that the operator $W[\partial\Omega ,{\mathbf{a}},S_{{\mathbf{a}}},\cdot]$
is bounded from $C^{0,\alpha}(\partial\Omega)$ to $C^{1,\alpha}(\partial\Omega)$ by the work of
Mikhlin and Pr\"{o}ssdorf \cite{MikPr86}.

In case  $n=3$, $m\geq 2$ and $\Omega$ is of class $C^{m,\alpha}$  and if  $P[{\mathbf{a}},D]$ is the  Helmholtz operator, Kirsch~\cite{Ki89} has proved that  the operator  $W_\Omega[{\mathbf{a}},S_{{\mathbf{a}}},\cdot]$  is  bounded from $C^{m-1,\alpha}(\partial\Omega)$ to $C^{m,\alpha}(\partial\Omega)$ and that accordingly it
 is  compact in $C^{m,\alpha}(\partial\Omega)$. 
 
 Then Heinemann~\cite{He92} has developed the ideas of  von Wahl in the frame of Schauder spaces and has proved that
 if $\Omega$ is of class $C^{m+5}$ and if 
 $S_{{\mathbf{a}}}$ is the fundamental solution of the Laplace operator,  then  the double layer improves the regularity of one unit on the boundary, \textit{i.e.}, 
 $W_\Omega[{\mathbf{a}},S_{{\mathbf{a}}},\cdot]$  is linear and continuous from $C^{m,\alpha}(\partial\Omega)$ to $C^{m+1,\alpha}(\partial\Omega)$.

Mitrea~\cite{Mit14}  has proved that the double layer of second order equations and systems is compact in  $C^{0,\beta}(\partial\Omega)$   for $\beta\in]0, \alpha[$ and bounded  in $C^{0,\alpha}(\partial\Omega)$
  under the assumption that $\Omega$ is of class $C^{1,\alpha}$. Then by exploiting a formula for the tangential derivatives such results have been extended to  compactness  and boundedness results in  $C^{1,\beta}(\partial\Omega)$ and $C^{1,\alpha}(\partial\Omega)$, respectively.
  
  In \cite{DoLa17}, we have proved that if $m\geq 1$, $\beta\in]0,\alpha]$, $\alpha\in]0,1[$, then 
  $W_\Omega[{\mathbf{a}},S_{{\mathbf{a}}},\cdot]$  is linear and continuous from $C^{m,\beta}(\partial\Omega)$ to $C^{m,\alpha}(\partial\Omega)$ and a related result if we chose $\beta=0$.
  
  In this paper we plan to consider the case in which  $\Omega$ is of class 
   $C^{1,1}$  and show that if the  maximal function of the tangential gradient with respect to the first variable of the kernel of the double layer potential is bounded,	   then $W_\Omega[{\mathbf{a}},S_{{\mathbf{a}}},\cdot]$  is linear and continuous from $C^{0,\beta}(\partial\Omega)$ to $C^{1,\beta}(\partial\Omega)$ for $\beta\in]0,1[$ and is linear and continuous from $C^{0,1}(\partial\Omega)$ to the generalized Schauder space
$C^{1,\omega_{1}}(\partial\Omega)$ of functions with $1$-st order tangential  derivatives which satisfy a generalized $\omega_{1}$-H\"{o}lder condition with
\[
\omega_{1}(r)\sim r^{1}|\ln r| 
\qquad{\mathrm{as}}\ r\to 0\,,
\]
see Theorem \ref{thm:dllreggen}. Our proofs are based on   Theorem of \cite[Thm.~6.6]{La22b} on integral operators, that   
we report here in the case in which the domain of integration is a compact differentiable manifold, see Theorem \ref{thm:iokreg}.  Theorem \ref{thm:iokreg} requires that we can estimate the maximal function associated to the tangential gradient of the kernel of the double layer potential with respect to its first variable and that the same tangential gradient belongs to a certain class of kernels. 
 Then we prove the membership in the class of kernels by exploiting the imbedding and multiplication properties that we have highlighted and proved in \cite{La22b} and that we  report here in the special case we need, see Section \ref{sec:kecla}. Here we note that the properties of Section \ref{sec:kecla} actually simplify a  proof that would be otherwise long to explain.
 
 \section{Notation}\label{sec:notation}
 Let $M_n({\mathbb{R}})$ denote the set of $n\times n$ matrices with real entries. 
     $\delta_{l,j}$ denotes the Kronecker  symbol. Namely,  $\delta_{l,j}=1$ if $l=j$, $\delta_{l,j}=0$ if $l\neq j$, with $l,j\in {\mathbb{N}}$. $|A|$ denotes the operator norm of a matrix $A$, 
       $A^{t}$ denotes the transpose matrix of $A$.  
 We set
\begin{equation}\label{eq:balls}
{\mathbb{B}}_n(\xi,r)\equiv \left\{\eta\in {\mathbb{R}}^n:\,  |\xi-\eta|<r\right\}\,, 
\end{equation}
for all $(\xi,r)\in {\mathbb{R}}^n\times ]0,+\infty[$. If ${\mathbb{D}}$ is a subset of $ {\mathbb{R}}^n$, 
then we set
\[
B({\mathbb{D}})\equiv\left\{
f\in {\mathbb{C}}^{\mathbb{D}}:\,f\ \text{is\ bounded}
\right\}
\,,\quad
\|f\|_{B({\mathbb{D}})}\equiv\sup_{\mathbb{D}}|f|\qquad\forall f\in B({\mathbb{D}})\,.
\]
Then $C^0({\mathbb{D}})$ denotes the set of continuous functions from ${\mathbb{D}}$ to ${\mathbb{C}}$ and we introduce the subspace
$
C^0_b({\mathbb{D}})\equiv C^0({\mathbb{D}})\cap B({\mathbb{D}})
$
of $B({\mathbb{D}})$.  Let $\omega$ be a function from $[0,+\infty[$ to itself such that
\begin{eqnarray}
\nonumber
&&\qquad\qquad\omega(0)=0,\qquad \omega(r)>0\qquad\forall r\in]0,+\infty[\,,
\\
\label{om}
&&\qquad\qquad\omega\ {\text{is\   increasing,}}\ \lim_{r\to 0^{+}}\omega(r)=0\,,
\\
\nonumber
&&\qquad\qquad{\text{and}}\ \sup_{(a,t)\in[1,+\infty[\times]0,+\infty[}
\frac{\omega(at)}{a\omega(t)}<+\infty\,.
\end{eqnarray}
If $f$ is a function from a subset ${\mathbb{D}}$ of ${\mathbb{R}}^n$   to ${\mathbb{C}}$,  then we denote by   $|f:{\mathbb{D}}|_{\omega (\cdot)}$  the $\omega(\cdot)$-H\"older constant  of $f$, which is delivered by the formula   
\[
|f:{\mathbb{D}}|_{\omega (\cdot)
}
\equiv
\sup\left\{
\frac{|f( x )-f( y)|}{\omega(| x- y|)
}: x, y\in {\mathbb{D}} ,  x\neq
 y\right\}\,.
\]        
If $|f:{\mathbb{D}}|_{\omega(\cdot)}<\infty$, we say that $f$ is $\omega(\cdot)$-H\"{o}lder continuous. Sometimes, we simply write $|f|_{\omega(\cdot)}$  
instead of $|f:{\mathbb{D}}|_{\omega(\cdot)}$. The
subset of $C^{0}({\mathbb{D}} ) $  whose
functions  are
$\omega(\cdot)$-H\"{o}lder continuous    is denoted  by  $C^{0,\omega(\cdot)} ({\mathbb{D}})$
and $|f:{\mathbb{D}}|_{\omega(\cdot)}$ is a semi-norm on $C^{0,\omega(\cdot)} ({\mathbb{D}})$.  
Then we consider the space  $C^{0,\omega(\cdot)}_{b}({\mathbb{D}} ) \equiv C^{0,\omega(\cdot)} ({\mathbb{D}} )\cap B({\mathbb{D}} ) $ with the norm \[
\|f\|_{ C^{0,\omega(\cdot)}_{b}({\mathbb{D}} ) }\equiv \sup_{x\in {\mathbb{D}} }|f(x)|+|f|_{\omega(\cdot)}\qquad\forall f\in C^{0,\omega(\cdot)}_{b}({\mathbb{D}} )\,.
\] 
\begin{rem}
\label{rem:om4}
Let $\omega$ be as in (\ref{om}). 
Let ${\mathbb{D}}$ be a   subset of ${\mathbb{R}}^{n}$. Let $f$ be a bounded function from $ {\mathbb{D}}$ to ${\mathbb{C}}$, $a\in]0,+\infty[$.  Then,
\[
\label{rem:om5}
\sup_{x,y\in {\mathbb{D}},\ |x-y|\geq a}\frac{|f(x)-f(y)|}{\omega(|x-y|)}
\leq \frac{2}{\omega(a)} \sup_{{\mathbb{D}}}|f|\,.
\]
\end{rem}
In the case in which $\omega(\cdot)$ is the function 
$r^{\alpha}$ for some fixed $\alpha\in]0,1]$, a so-called H\"{o}lder exponent, we simply write $|\cdot:{\mathbb{D}}|_{\alpha}$ instead of
$|\cdot:{\mathbb{D}}|_{r^{\alpha}}$, $C^{0,\alpha} ({\mathbb{D}})$ instead of $C^{0,r^{\alpha}} ({\mathbb{D}})$, $C^{0,\alpha}_{b}({\mathbb{D}})$ instead of $C^{0,r^{\alpha}}_{b} ({\mathbb{D}})$, and we say that $f$ is $\alpha$-H\"{o}lder continuous provided that 
$|f:{\mathbb{D}}|_{\alpha}<\infty$.
 \section{Special classes of potential type kernels in ${\mathbb{R}}^n$}\label{sec:kecla}
 In this section we collect some basic properties of the classes of kernel that we need. For the proofs, we refer to \cite[\S 3]{La22b}. If $X$ and $Y$ are subsets of ${\mathbb{R}}^n$, then we denote by ${\mathbb{D}}_{X\times Y}$ the diagonal  of $X\times Y$, i.e., we set
 \begin{equation}\label{diagonal}
{\mathbb{D}}_{X\times Y}\equiv\left\{
(x,y)\in X\times Y:\,x=y
\right\} 
\end{equation}
and if $X=Y$, then    we denote by ${\mathbb{D}}_{X}$ the diagonal  of $X\times X$, i.e., we set
\[
{\mathbb{D}}_X\equiv {\mathbb{D}}_{X\times X}\,.
\]
An off-diagonal function  in $X\times Y$ is a function from $(X\times Y)\setminus {\mathbb{D}}_{X\times Y}$ to ${\mathbb{C}}$. We now wish to consider a specific class of off-diagonal kernels.\begin{defn}
 Let $X$ and $Y$ be subsets of ${\mathbb{R}}^n$. Let $s\in {\mathbb{R}}$. We denote by ${\mathcal{K}}_{s,X\times Y}$ (or more simply by ${\mathcal{K}}_s$), the set of continuous functions $K$ from $(X\times Y)\setminus {\mathbb{D}}_{ X\times Y }$ to ${\mathbb{C}}$ such that
\[
 \|K\|_{ {\mathcal{K}}_{s,X\times Y} }\equiv \sup_{(x,y)\in  (X\times Y)\setminus {\mathbb{D}}_{ X\times Y }  }|K(x,y)|\,|x-y|^s<+\infty\,.
\]
The elements of $ {\mathcal{K}}_{s,X\times Y}$ are said to be kernels of potential type $s$ in $X\times Y$. 
\end{defn}
We plan to consider `potential type' kernels as in the following definition
(see also paper \cite{DoLa17} with Dondi,  where such classes have been introduced in a form that generalizes those of Gegelia \cite{Ge67}, \cite[Chap.~IV]{KuGeBaBu79} and Giraud \cite{Gi34}).	 
\begin{defn}\label{defn:ksss}
 Let $X$, $Y\subseteq {\mathbb{R}}^n$. Let $s_1$, $s_2$, $s_3\in {\mathbb{R}}$. We denote by ${\mathcal{K}}_{s_1, s_2, s_3} (X\times Y)$ the set of continuous functions $K$ from $(X\times Y)\setminus {\mathbb{D}}_{X\times Y}$ to ${\mathbb{C}}$ such that
 \begin{eqnarray*}
\lefteqn{
\|K\|_{  {\mathcal{K}}_{ s_1, s_2, s_3  }(X\times Y)  }
\equiv
\sup\biggl\{\biggr.
|x-y|^{ s_{1} }|K(x,y)|:\,(x,y)\in X\times Y, x\neq y
\biggl.\biggr\}
}
\\ \nonumber
&&\qquad\qquad\qquad
+\sup\biggl\{\biggr.
\frac{|x'-y|^{s_{2}}}{|x'-x''|^{s_{3}}}
|  K(x',y)- K(x'',y)  |:\,
\\ \nonumber
&&\qquad\qquad\qquad 
x',x''\in X, x'\neq x'', y\in Y\setminus{\mathbb{B}}_{n}(x',2|x'-x''|)
\biggl.\biggr\}<+\infty\,.
\end{eqnarray*}
\end{defn}
 One can easily verify that $({\mathcal{K}}_{ s_{1},s_{2},s_{3}   }(X\times Y),\|\cdot\|_{  {\mathcal{K}}_{s_{1},s_{2},s_{3}   }(X\times Y)  })$ is a normed space.  By our definition, if $s_1$, $s_2$, $s_3\in {\mathbb{R}}$, we have
\[
{\mathcal{K}}_{s_{1},s_{2},s_{3}   }(X\times Y) \subseteq {\mathcal{K}}_{s_{1}, X\times Y} 
\]
and
\[
\|K\|_{{\mathcal{K}}_{s_{1}, X\times Y} }\leq \|K\|_{ {\mathcal{K}}_{s_{1},s_{2},s_{3}   }(X\times Y) }
\qquad\forall K\in {\mathcal{K}}_{s_{1},s_{2},s_{3}   }(X\times Y) \,.
\]
We note that if we choose $s_2=s_1+s_3$ we have a so-called class of standard kernels.   Then we have the following elementary known embedding lemma (cf.~\textit{e.g.}, \cite[\S 3]{La22b}).
\begin{lem}\label{lem:kelem}
 Let $X$, $Y\subseteq {\mathbb{R}}^n$. Let $s_1$, $s_2$, $s_3\in {\mathbb{R}}$. If $a\in]0,+\infty[$, then $ {\mathcal{K}}_{s_{1},s_{2},s_{3}   }(X\times Y)$ is continuously embedded into  ${\mathcal{K}}_{s_{1},s_{2}-a,s_{3}-a   }(X\times Y)$.
\end{lem}
Next we introduce the following known elementary lemma, which we exploit later and which can be proved by the triangular inequality.
\begin{lem}\label{lem:rec}
 \[
 \frac{1}{2 }|x'-y|\leq |x''-y|\leq 2|x'-y|\,,
 \]
 for all  $x',x''\in {\mathbb{R}}^{n}$, $x'\neq x''$, $y\in {\mathbb{R}}^{n}
\setminus {\mathbb{B}}_{n}(x',2|x'-x''|)$.
\end{lem}
Next we state the following two product rule statements (cf.~\cite[\S 3]{La22b}).
\begin{thm}\label{thm:kerpro}
 Let $X$, $Y\subseteq {\mathbb{R}}^n$. Let $s_1$, $s_2$, $s_3$, $t_1$, $t_2$, $t_3\in {\mathbb{R}}$. 
\begin{enumerate}
\item[(i)] If $K_1\in  {\mathcal{K}}_{s_{1},s_{2},s_{3}   }(X\times Y)$ and $K_2\in  {\mathcal{K}}_{t_{1},t_{2},t_{3}   }(X\times Y)$, then the following inequality holds
\begin{eqnarray*}
\lefteqn{
|  K_1(x',y)K_2(x',y)- K_1(x'',y)K_2(x'',y)  |
}
\\ \nonumber
&&\qquad\qquad\qquad\qquad \qquad 
\leq\|K_1\|_{ {\mathcal{K}}_{s_{1},s_{2},s_{3}   }(X\times Y)}
\|K_2\|_{ {\mathcal{K}}_{t_{1},t_{2},t_{3}   }(X\times Y)}
\\ \nonumber
&&\qquad\qquad\qquad\qquad \qquad   
\quad
\times\left(
\frac{|x'-x''|^{s_3}}{|x'-y|^{s_2+t_1}}+ \frac{2^{|s_1|}|x'-x''|^{t_3}}{|x'-y|^{t_2+s_1}}
\right)
\end{eqnarray*}
for all $x',x''\in X$, $x'\neq x''$, $y\in Y\setminus{\mathbb{B}}_{n}(x',2|x'-x''|)$. 
\item[(ii)] The pointwise product is bilinear and continuous from
\[
{\mathcal{K}}_{s_{1},s_1+s_3,s_{3}   }(X\times Y)\times  {\mathcal{K}}_{t_{1},t_{1}+s_3,s_{3}   }(X\times Y)
\quad\text{to}\quad
{\mathcal{K}}_{s_1+t_{1},s_{1}+s_3+t_1,s_{3}   }(X\times Y)\,.
\] 
\end{enumerate}
\end{thm}

\begin{prop}\label{prop:prkerho}
 Let $X$, $Y\subseteq {\mathbb{R}}^n$. Let $s_1$, $s_2$, $s_3 \in {\mathbb{R}}$,  
 $\alpha\in]0,1]$. Then the following statements hold.
\begin{enumerate}
\item[(i)] If $K\in {\mathcal{K}}_{s_{1},s_2,s_{3}   }(X\times Y)$ and $f\in C^{0,\alpha}_b(X)$, then
\[
|K(x,y)f(x)|\,|x-y|^{s_1}\leq \|K\|_{ {\mathcal{K}}_{s_{1},X\times Y  } }\sup_X|f|
\qquad\forall (x,y)\in X\times Y\setminus{\mathbb{D}}_{X\times Y}\,.
\]
and
\begin{eqnarray*}
\lefteqn{
|K(x',y)f(x')-K(x'',y)f(x'')|
}
\\ \nonumber
&&\qquad
\leq \|K\|_{	 {\mathcal{K}}_{s_{1},s_2,s_{3}   }(X\times Y) }\|f\|_{	C^{0,\alpha}_b(X)	}
\left\{
\frac{|x'-x''|^{s_3}}{|x'-y|^{s_2}}+2^{|s_1|}\frac{|x'-x''|^{\alpha}}{|x'-y|^{s_1}} 
\right\}
\end{eqnarray*}
for all $x',x''\in X$, $x'\neq x''$, $y\in Y\setminus{\mathbb{B}}_{n}(x',2|x'-x''|)$. 
\item[(ii)] If  $s_2\geq s_1$ and $X$ and $Y$ are both bounded, then the map from 
\[
{\mathcal{K}}_{s_{1},s_2,s_{3}   }(X\times Y)\times C^{0,s_3}_b(X)\quad\text{to}\quad{\mathcal{K}}_{s_{1},s_2,s_{3}   }(X\times Y)
\]
 that takes the pair $(K,f)$ to the kernel $K(x,y)f(x)$ of the variable $(x,y)\in (X\times Y)\setminus {\mathbb{D}}_{X\times Y}$ is bilinear and continuous. 
 \item[(iii)] The map from 
\[
{\mathcal{K}}_{s_{1},s_2,s_{3}   }(X\times Y)\times C^{0}_b(Y)\quad\text{to}\quad{\mathcal{K}}_{s_{1},s_2,s_{3}   }(X\times Y)
\]
 that takes the pair $(K,f)$ to the kernel $K(x,y)f(y)$ of the variable $(x,y)\in (X\times Y)\setminus {\mathbb{D}}_{X\times Y}$ is bilinear and continuous. 
\end{enumerate}
\end{prop}
Next we have the following imbedding statement that holds for bounded sets (cf.~\cite[\S 3]{La22b}).
\begin{prop}\label{prop:kerem}
Let $X$, $Y$ be bounded subsets of ${\mathbb{R}}^n$. Let $s_1$, $s_2$, $s_3$, $t_1$, $t_2$, $t_3\in {\mathbb{R}}$. Then the following statements hold.
\begin{enumerate}
\item[(i)] If $t_1\geq s_1$ then $ {\mathcal{K}}_{s_1,X\times Y}$ is continuously embedded into  $ {\mathcal{K}}_{t_1,X\times Y}$.
\item[(ii)] If $t_1\geq s_1$, $t_3\leq s_3$ and $(t_2-t_3)\geq(s_2-s_3)$, then ${\mathcal{K}}_{s_{1},s_{2},s_{3}   }(X\times Y)$ is continuously embedded into ${\mathcal{K}}_{t_{1},t_{2},t_{3}   }(X\times Y)$.
\item[(iii)] If $t_1\geq s_1$,  $t_3\leq s_3$, then ${\mathcal{K}}_{s_{1},s_{1}+s_3,s_{3}   }(X\times Y)$ is continuously embedded into the space ${\mathcal{K}}_{t_{1},t_{1}+t_3,t_{3}   }(X\times Y)$.
\end{enumerate}
\end{prop}
We now show that we can associate a potential type kernel to all H\"{o}lder continuous functions
(cf.~\cite[\S 3]{La22b}).
\begin{lem}\label{lem:hoker}
 Let $X$, $Y$ be subsets of ${\mathbb{R}}^n$. Let $\alpha\in]0,1]$. Then the following statements hold.
 \begin{enumerate}
\item[(i)] If $\mu\in C^{0,\alpha}(X\cup Y)$, then the map $\Xi[\mu]$  defined by
\begin{equation}\label{lem:hoker1}
\Xi[\mu](x,y)\equiv \mu(x)-\mu(y)\qquad\forall (x,y)\in (X\times Y)\setminus {\mathbb{D}}_{X\times Y}
\end{equation}
 belongs to ${\mathcal{K}}_{-\alpha,0,\alpha  }(X\times Y)$.
\item[(ii)] The operator $\Xi$ from $C^{0,\alpha}(X\cup Y)$ to ${\mathcal{K}}_{-\alpha,0,\alpha  }(X\times Y)$ that takes $\mu$ to $\Xi[\mu]$ is linear and continuous.
\end{enumerate}
 \end{lem}
In order to introduce a result  of \cite[Thm.~6.6]{La22b}, we need to introduce a further norm for kernels
 in the case in which  $ Y$ is a compact manifold of class $C^1$ that is imbedded in $ M=	  {\mathbb{R}}^n$ and $X= Y$.
\begin{defn}
Let $ Y$ be a compact manifold of class $C^1$ that is imbedded in ${\mathbb{R}}^n$.
    Let $s_1$, $s_2$, $s_3\in {\mathbb{R}}$. We set
\begin{eqnarray*}
\lefteqn{
 {\mathcal{K}}_{ s_1, s_2, s_3  }^\sharp( Y	\times  Y	)
  \equiv 
  \biggl\{\biggr.
K\in  {\mathcal{K}}_{ s_1, s_2, s_3  }( Y	 \times  Y	):\,
}
\\ \nonumber
&&\ \ 
\sup_{x\in  Y	}\sup_{r\in ]0,+\infty[}
\left|
\int_{ Y	\setminus {\mathbb{B}}_n(x,r)}K(x,y)\,d\nu(y)
\right|<+\infty
 \biggl.\biggr\}  
\end{eqnarray*}
and
\begin{eqnarray*}
\lefteqn{
\|K\|_{{\mathcal{K}}_{ s_1, s_2, s_3  }^\sharp( Y	\times Y	)}
\equiv
\|K\|_{{\mathcal{K}}_{ s_1, s_2, s_3  }( Y	\times Y	)}
}
\\ \nonumber
&&+
\sup_{x\in Y}\sup_{r\in ]0,+\infty[}
\left|
\int_{Y\setminus {\mathbb{B}}_n(x,r)}K(x,y)\,d\nu(y)
\right|\quad\forall 
 K\in {\mathcal{K}}_{ s_1, s_2, s_3  }^\sharp( Y	\times  Y	 )\,.
 \end{eqnarray*}
\end{defn}
Clearly,  $({\mathcal{K}}^\sharp_{ s_{1},s_{2},s_{3}   }( Y\times  Y),\|\cdot\|_{  {\mathcal{K}}^\sharp_{s_{1},s_{2},s_{3}   }(Y	\times Y)  })$ is a normed space. By definition, ${\mathcal{K}}^\sharp_{ s_{1},s_{2},s_{3}   }( Y	\times Y	)$ is continuously embedded into  ${\mathcal{K}}_{ s_{1},s_{2},s_{3}   }( Y	\times Y	)$. Next we introduce a function that we need for a generalized H\"{o}lder norm. For each $\theta\in]0,1]$, we define the function $\omega_{\theta}(\cdot)$ from $[0,+\infty[$ to itself by setting
\[
\omega_{\theta}(r)\equiv
\left\{
\begin{array}{ll}
0 &r=0\,,
\\
r^{\theta}|\ln r | &r\in]0,r_{\theta}]\,,
\\
r_{\theta}^{\theta}|\ln r_{\theta} | & r\in ]r_{\theta},+\infty[\,,
\end{array}
\right.
\]
where
$
r_{\theta}\equiv e^{-1/\theta}
$ for all $\theta\in ]0,1]$. Obviously, $\omega_{\theta}(\cdot) $ is concave and satisfies   condition (\ref{om}).
 We also note that if ${\mathbb{D}}\subseteq {\mathbb{R}}^n$, then the continuous embedding
\[
C^{0, \theta }_b({\mathbb{D}})\subseteq 
C^{0,\omega_\theta(\cdot)}_b({\mathbb{D}})\subseteq 
C^{0,\theta'}_b({\mathbb{D}})
\]
holds for all $\theta'\in ]0,\theta[$. Here the subscript $b$ denotes that we are considering the intersection of a (generalized) H\"older space with the space $B({\mathbb{D}})$ of the bounded functions in ${\mathbb{D}}$. Then we introduce the following result of \cite[Thm.~6.3
]{La22b}.
 \begin{thm}\label{thm:iokreg}
 Let $ Y	 $ be a compact manifold of class $C^1$ that is imbedded in ${\mathbb{R}}^n$.  Let $s_1\in [0,(n-1)[$. 
 Let $\beta\in ]0,1]$, $t_1\in[\beta,(n-1)+\beta[$, $t_2\in [ \beta,+\infty[$, $t_3\in]0,1]$. Let the kernel $K\in {\mathcal{K}}_{s_1,s_1+1,1}( Y	 \times  Y	)$ satisfy the following assumption 
 \[
 K(\cdot,y)\in C^1( Y	 \setminus\{y\})  \quad\forall y\in  Y	\,.
 \]
Then the following statements hold.
\begin{enumerate}
\item[(i)] If $t_1<(n-1)$ and ${\mathrm{grad}}_{ Y	,x}K(\cdot,\cdot)\in\left(
 {\mathcal{K}}_{t_1,t_2,t_3}( Y	\times  Y	)
 \right)^n$, then the following statements hold.

\begin{enumerate}
\item[(a)] If $t_2-\beta>(n-1)$, $t_2<(n-1)+\beta+t_3$ and 
\[
 \int_Y	K(\cdot,y)\,  	d\sigma_y
\in C^{1,\min\{ \beta, (n-1)+t_3+\beta-t_2\} }( Y	)\,,
\]
then the map from $C^{0,\beta} ( Y	)$ to $C^{1,\min\{ \beta, (n-1)+t_3+\beta-t_2\}} ( Y	)$ that takes $\mu$ to the function  $\int_YK(\cdot,y)\mu(y)\,d\sigma_y$ is linear and continuous. 
\item[(aa)] If $t_2-\beta=(n-1)$  and 
\[
\int_Y	K(\cdot,y)\, 	d\sigma_y
\in C^{1,\max\{ r^\beta, \omega_{t_3} (\cdot)\} }( Y	)\,,
\] 
then the map from $C^{0,\beta}( Y	)$ to $C^{1,\max\{ r^\beta, \omega_{t_3} (\cdot)\} }( Y	)
$ that takes $\mu$ to the function $\int_Y	K(\cdot,y)\mu(y)\,d\sigma_y$ is linear and continuous. 
\end{enumerate}

\item[(ii)] If $t_1=(n-1)$ and ${\mathrm{grad}}_{Y	,x}K(\cdot,\cdot)\in\left(
 {\mathcal{K}}^\sharp_{t_1,t_2,t_3}(Y	\times Y	)
 \right)^n$, then the following statements hold.

\begin{enumerate}
\item[(b)] If $t_2-\beta>(n-1)$, $t_2<(n-1)+\beta+t_3$ and 
\[
 \int_Y	K(\cdot,y)\, d\sigma_y 
 \in C^{1,\min\{  \beta, (n-1)+t_3+\beta-t_2\} }( Y	)\,,
\]
then the map from $C^{0,\beta}( Y	)$ to $C^{1,\min\{  \beta, (n-1)+t_3+\beta-t_2\}}_b( Y)$ that takes $\mu$ to the function  $\int_Y	K(\cdot,y)\mu(y)\,d\sigma_y$ is linear and continuous. 
\item[(bb)] If $t_2-\beta=(n-1)$  and 
\[
\int_Y	K(\cdot,y)\,d\sigma_y
\in C^{1,\max\{r^\beta, \omega_{t_3} (\cdot)\} }(Y)\,,
\] 
then the map from $C^{0,\beta}(Y)$ to $C^{1,\max\{ r^\beta, \omega_{t_3} (\cdot)\} }(Y)
$ that takes $\mu$ to the function $\int_Y	K(\cdot,y)\mu(y)\,d\sigma_y$ is linear and continuous. 
\end{enumerate}

\item[(iii)] If $t_1>(n-1)$ and ${\mathrm{grad}}_{Y	,x}K(\cdot,\cdot)\in\left(
 {\mathcal{K}}_{t_1,t_2,t_3}( Y	\times  Y	)
 \right)^n$, then the following statements hold.

\begin{enumerate}
\item[(c)] If $t_2-\beta>(n-1)$, $t_2<(n-1)+\beta+t_3$ and 
\[
 \int_YK(\cdot,y)\, d\sigma_y 
 \in C^{1,\min\{\beta,(n-1)+ \beta-t_1, (n-1)+t_3+\beta-t_2\} }( Y)\,,
\]
then the map from $C^{0,\beta} (Y)$ to $C^{1,\min\{\beta,(n-1)+ \beta-t_1, (n-1)+t_3+\beta-t_2\}}(Y)$ that takes $\mu$ to the function  $\int_YK(\cdot,y)\mu(y)\,d\sigma_y$ is linear and continuous. 
\item[(cc)] If $t_2-\beta=(n-1)$  and 
\[
\int_YK(\cdot,y)\, d\sigma_y 
\in C^{1,\max\{ r^\beta,r^{(n-1)+ \beta-t_1}, \omega_{t_3} (\cdot)\} }( Y)\,,
\] 
then the map from $C^{0,\beta}(Y)$ to $C^{1,\max\{r^\beta, r^{(n-1)+ \beta-t_1}, \omega_{t_3} (\cdot)\} }( Y	)
$ that takes $\mu$ to the function $\int_Y	K(\cdot,y)\mu(y)\,d\sigma_y$ is linear and continuous. 
\end{enumerate}
\end{enumerate}
\end{thm}
 We also need to consider convolution kernels, thus we introduce the following notation.
 If $n\in {\mathbb{N}}\setminus\{0\}$, $m\in {\mathbb{N}}$, $h\in {\mathbb{R}}$, $\alpha\in ]0,1]$, then we set 
\begin{equation}\label{volume.klhomog}
{\mathcal{K}}^{m,\alpha}_h \equiv\biggl\{
k\in C^{m,\alpha}_{ {\mathrm{loc}}}({\mathbb{R}}^n\setminus\{0\}):\, k\ {\text{is\ positively\ homogeneous\ of \ degree}}\ h
\biggr\}\,,
\end{equation}
where $C^{m,\alpha}_{ {\mathrm{loc}}}({\mathbb{R}}^n\setminus\{0\})$ denotes the set of functions of 
$C^{m}({\mathbb{R}}^n\setminus\{0\})$ whose restriction to $\overline{\Omega}$ is of class $C^{m,\alpha}(\overline{\Omega})$ for all bounded open subsets $\Omega$ of ${\mathbb{R}}^n$ such that
$\overline{\Omega}\subseteq {\mathbb{R}}^n\setminus\{0\}$
and we set
\[
\|k\|_{ {\mathcal{K}}^{m,\alpha}_h}\equiv \|k\|_{C^{m,\alpha}(\partial{\mathbb{B}}_n(0,1))}\qquad\forall k\in {\mathcal{K}}^{m,\alpha}_h\,.
\]
We can easily verify that $ \left({\mathcal{K}}^{m,\alpha}_h , \|\cdot\|_{ {\mathcal{K}}^{m,\alpha}_h}\right)$ is a Banach space. We also mention the following variant of a well known statement.
 \begin{lem}\label{lem:knnnone}
 Let $n\in {\mathbb{N}}\setminus\{0\}$, $h\in [0,+\infty[$. If $k\in C^{0,1}_{ {\mathrm{loc}} }({\mathbb{R}}^n\setminus\{0\})$ is positively homogeneous of degree $-h$, then $k(x-y)\in {\mathcal{K}}_{h,h+1,1}({\mathbb{R}}^n\times {\mathbb{R}}^n)$.
 Moreover, the map from ${\mathcal{K}}^{0,1}_{-h}$   to ${\mathcal{K}}_{h,h+1,1}({\mathbb{R}}^n\times {\mathbb{R}}^n)$ which takes $k$ to $k(x-y)$ is linear and continuous (see (\ref{volume.klhomog}) for the definition of $
{\mathcal{K}}^{0,1}_{-h}$).
\end{lem}
{\bf Proof.} Since $k$ is positively homogeneous  of degree $-h$, we have
 \[
 |k(x-y)|\leq (\sup_{\partial{\mathbb{B}}_n(0,1)}|k|)
 |x-y|^{-h}\qquad\forall (x,y)\in ({\mathbb{R}}^n\times {\mathbb{R}}^n)\setminus {\mathbb{D}}_{{\mathbb{R}}^n\times {\mathbb{R}}^n}\,.
 \]
 Since $k$ is positively homogeneous of degree $-(n-1)$, the inequality of Cialdea \cite[VIII, p.~47]{Ci95} (see also 
\cite[Lem.~4.14]{DaLaMu21}  with   $\alpha=1$) implies that if $x'$, $x''\in {\mathbb{R}}^n$, $x'\neq x''$, $y\in  {\mathbb{R}}^n\setminus{\mathbb{B}}_n(x',2|x'-x''|)$, then
\begin{eqnarray*}
\lefteqn{
|k(x'-y)-k(x''-y)|
}
\\ \nonumber
&&\qquad
\quad
\leq  ( 2^1+ 2h)
 \max\{\sup_{ \partial{\mathbb{B}}_n(0,1) }|k|, |k:\,\partial{\mathbb{B}}_n(0,1)|_1\}
 \\ \nonumber
&&\qquad
\quad\quad
 \times |(x'-y)-(x''-y)| (\min\{|(x'-y)|,|(x''-y)|\})^{-h-1}\,.
\end{eqnarray*}
Then Lemma \ref{lem:rec} implies that
$
|x''-y|\geq \frac{1}{2}|x'-y|$, 
and thus we have
\begin{eqnarray*}
\lefteqn{
|k(x'-y)-k(x''-y)|
}
\\ \nonumber
&&\qquad
\leq (2+2h)\max\{\sup_{ \partial{\mathbb{B}}_n(0,1) }|k|, |k:\,\partial{\mathbb{B}}_n(0,1)|_1\}
\frac{|x'-x''|}{|x'-y|^{h+1}}2^{h+1}
\end{eqnarray*}
and the proof is complete.\hfill  $\Box$ 

\vspace{\baselineskip}

If $X$ and $Y$ are subsets of ${\mathbb{R}}^n$, then the restriction operator  
\[
\text{from}\ {\mathcal{K}}_{h,h+1,1}({\mathbb{R}}^n\times {\mathbb{R}}^n)\ \text{to}\
{\mathcal{K}}_{h,h+1,1}(X\times Y)
\]
is linear and continuous. Thus Lemma \ref{lem:knnnone} implies that the map
\[
\text{from\ the subspace}\  {\mathcal{K}}^{0,1}_{-h}\ \text{ 
 of}\  C^{0,1}_{ {\mathrm{loc}} }({\mathbb{R}}^n\setminus\{0\})\ \text{ to}\   {\mathcal{K}}_{h,h+1,1}(X\times Y)\,,
 \]
  which takes $k$ to $k(x-y)$ is linear and continuous.

\begin{rem} 
 As Lemma \ref{lem:knnnone} shows the convolution kernels associated to positively homogeneous functions of negative degree are standard kernels. We note however that 
 there exist potential type kernels that belong to a class ${\mathcal{K}}_{s_1,s_2,s_3} (X\times Y)$ with   $s_2\neq s_1+s_3$. 
\end{rem}

\section{Technical preliminaries on the differential operator}\label{sec:techdiop}
 Let $\Omega$ be a bounded open subset of ${\mathbb{R}}^{n}$ of class $C^{2}$.  
 The kernel of the boundary integral operator corresponding to the  double layer potential 
  is the following
\begin{eqnarray}\label{eq:tgdlgen}
\lefteqn{
\overline{B^{*}_{\Omega,y}}\left(S_{{\mathbf{a}}}(x-y)\right)\equiv - \sum_{l,j=1}^{n} a_{jl}\nu_{l}(y)\frac{\partial S_{ {\mathbf{a}} } }{\partial x_{j}}(x-y) 
}
\\ \nonumber
&&\qquad\qquad\qquad
 - \sum_{l=1}^{n}\nu_{l}(y)a_{l}
S_{ {\mathbf{a}} }(x-y)\qquad\forall (x,y)\in (\partial\Omega)^2\setminus {\mathbb{D}}_{\partial\Omega} 
\end{eqnarray}
(cf.~(\ref{introd3})). In order to analyze the kernel of the double layer potential, we need some more information on the fundamental solution $S_{ {\mathbf{a}} } $. 
 To do so, we introduce the fundamental solution $S_{n}$ of the Laplace operator. Namely, we set
\[
S_{n}(x)\equiv
\left\{
\begin{array}{lll}
\frac{1}{s_{n}}\ln  |x| \qquad &   \forall x\in 
{\mathbb{R}}^{n}\setminus\{0\},\quad & {\mathrm{if}}\ n=2\,,
\\
\frac{1}{(2-n)s_{n}}|x|^{2-n}\qquad &   \forall x\in 
{\mathbb{R}}^{n}\setminus\{0\},\quad & {\mathrm{if}}\ n>2\,,
\end{array}
\right.
\]
where $s_{n}$ denotes the $(n-1)$ dimensional measure of 
$\partial{\mathbb{B}}_{n}(0,1)$ and	 
 we follow a formulation of Dalla Riva \cite[Thm.~5.2, 5.3]{Da13} and Dalla Riva, Morais and Musolino \cite[Thm.~5.5]{DaMoMu13}, that we state  as in paper  \cite[Cor.~4.2]{DoLa17} with Dondi (see also John~\cite{Jo55}, and Miranda~\cite{Mi65} for homogeneous operators, and Mitrea and Mitrea~\cite[p.~203]{MitMit13}).   
\begin{prop}
 \label{prop:ourfs} 
Let ${\mathbf{a}}$ be as in (\ref{introd0}), (\ref{ellip}), (\ref{symr}). 
Let $S_{ {\mathbf{a}} }$ be a fundamental solution of $P[{\mathbf{a}},D]$. 
Then there exist an invertible matrix $T\in M_{n}({\mathbb{R}})$ such that
\begin{equation}
\label{prop:ourfs0}
a^{(2)}=TT^{t}\,,
\end{equation}
 a real analytic function $A_{1}$ from $\partial{\mathbb{B}}_{n}(0,1)\times{\mathbb{R}}$ to ${\mathbb{C}}$ such that
 $A_{1}(\cdot,0)$ is odd,    $b_{0}\in {\mathbb{C}}$, a real analytic function $B_{1}$ from ${\mathbb{R}}^{n}$ to ${\mathbb{C}}$ such that $B_{1}(0)=0$, and a real analytic function $C $ from ${\mathbb{R}}^{n}$ to ${\mathbb{C}}$ such that
\begin{eqnarray}
\label{prop:ourfs1}
\lefteqn{
S_{ {\mathbf{a}} }(x)
= 
\frac{1}{\sqrt{\det a^{(2)} }}S_{n}(T^{-1}x)
}
\\ \nonumber
&&\qquad
+|x|^{3-n}A_{1}(\frac{x}{|x|},|x|)
 +(B_{1}(x)+b_{0}(1-\delta_{2,n}))\ln  |x|+C(x)\,,
\end{eqnarray}
for all $x\in {\mathbb{R}}^{n}\setminus\{0\}$,
 and such that both $b_{0}$ and $B_{1}$   equal zero
if $n$ is odd. Moreover, 
 \[
 \frac{1}{\sqrt{\det a^{(2)} }}S_{n}(T^{-1}x) 
 \]
is a fundamental solution for the principal part
  of $P[{\mathbf{a}},D]$.
\end{prop}
In particular for the statement that $A_{1}(\cdot,0)$ is odd, we refer to
Dalla Riva, Morais and Musolino \cite[Thm.~5.5, (32)]{DaMoMu13}, where $A_{1}(\cdot,0)$ coincides with ${\mathbf{f}}_1({\mathbf{a}},\cdot)$ in that paper. Here we note that a function $A$ from $(\partial{\mathbb{B}}_{n}(0,1))\times{\mathbb{R}}$ to ${\mathbb{C}}$ is said to be real analytic provided that it has a real analytic extension   to an open neighbourhood of $(\partial{\mathbb{B}}_{n}(0,1))\times{\mathbb{R}}$ in 
${\mathbb{R}}^{n+1}$. Then we have the following elementary lemma. 
\begin{lem}\label{lem:anexsph}
 Let $n\in {\mathbb{N}}\setminus\{0,1\}$. A function $A$ from   $(\partial{\mathbb{B}}_{n}(0,1))\times{\mathbb{R}}$ to ${\mathbb{C}}$ is  real analytic if and only if the function $\tilde{A}$ from $({\mathbb{R}}^n\setminus\{0\}) \times{\mathbb{R}}$ defined by 
\begin{equation}\label{lem:anexsph1}
\tilde{A}(x,r)\equiv A(\frac{x}{|x|},r)\qquad\forall (x,r)\in ({\mathbb{R}}^n\setminus\{0\}) \times{\mathbb{R}}
\end{equation}
is real analytic.
 \end{lem}  
 {\bf Proof.} If $A$ is real  analytic then, it has a real analytic extension $A^\sharp$ to an open neighborhood $U$ of $(\partial{\mathbb{B}}_{n}(0,1))\times{\mathbb{R}}$ in 
${\mathbb{R}}^{n+1}$. Since the function $\frac{x}{|x|}$ is real analytic in $x\in {\mathbb{R}}^n\setminus\{0\}$, then the composition $\tilde{A}$
of $A^\sharp$ and of $(\frac{x}{|x|},r)$ is real analytic. 

Conversely, if $\tilde{A}$ is real analytic, we note that 
$\tilde{A}$ is an extension of $A$ to the open neighborhood $({\mathbb{R}}^n\setminus\{0\}) \times{\mathbb{R}}$ of $(\partial{\mathbb{B}}_{n}(0,1))\times{\mathbb{R}}$  in 
${\mathbb{R}}^{n+1}$ and that accordingly $A$ is real analytic. \hfill  $\Box$ 

\vspace{\baselineskip}

Then  one can prove the following formula for the gradient of the fundamental solution (see reference \cite[Lem.~4.3,  (4.8) and the following 2 lines]{DoLa17} with Dondi. Here one should remember that $A_1(\cdot,0)$ is odd and that $b_0=0$ if $n$ is odd). 
\begin{prop}
\label{prop:grafun}
 Let ${\mathbf{a}}$ be as in (\ref{introd0}), (\ref{ellip}), (\ref{symr}). Let $T\in M_{n}({\mathbb{R}})$  be as in (\ref{prop:ourfs0}). Let $S_{ {\mathbf{a}} }$ be a fundamental solution of $P[{\mathbf{a}},D]$. Let  $B_{1}$, $C$
 be as in Proposition \ref{prop:ourfs}. 
  Then there exists a real analytic function $A_{2}$ from $\partial{\mathbb{B}}_{n}(0,1)\times{\mathbb{R}}$ to ${\mathbb{C}}^{n}$ such that
\begin{eqnarray}
\label{grafun1}
\lefteqn{
DS_{ {\mathbf{a}} }(x)=\frac{1}{ s_{n}\sqrt{\det a^{(2)} } }
|T^{-1}x|^{-n}x^{t}(a^{(2)})^{-1}
}
\\ \nonumber
&& 
+|x|^{2-n}A_{2}(\frac{x}{|x|},|x|)+DB_{1}(x)\ln |x|+DC(x)
\quad\forall x\in {\mathbb{R}}^{n}\setminus\{0\}\,.
\end{eqnarray}
Moreover,   $A_2(\cdot,0)$ is even.
\end{prop}
Then one can prove the following formula 
for the kernel of the double layer potential
\begin{eqnarray}
\label{eq:boest3}
\lefteqn{
\overline{B^{*}_{\Omega,y}}\left(S_{{\mathbf{a}}}(x-y)\right)
= -DS_{{\mathbf{a}}}(x-y)    a^{(2)}  \nu  (y)
- 
\nu^{t} (y)
a^{(1)}S_{{\mathbf{a}}}(x-y)
}
\\ \nonumber
&&\qquad\qquad
= 
-\frac{1}{s_{n}\sqrt{\det a^{(2)} }}|T^{-1}(x-y)|^{-n} (x-y)^{t}\nu (y)
\\ \nonumber
&&\qquad\qquad\quad
-
|x-y|^{2-n}
A_{2}(\frac{x-y}{|x-y|},|x-y|)a^{(2)}  \nu  (y)
\\ \nonumber
&&\qquad\qquad\quad
-DB_{1}(x-y)a^{(2)}  \nu  (y)\ln |x-y|
-   DC(x-y)a^{(2)}  \nu  (y)
\\ \nonumber
&&\qquad\qquad\quad
-\nu^{t}(y) a^{(1)}S_{{\mathbf{a}}}(x-y)
\qquad\forall x,y\in\partial\Omega, x\neq y\,.
\end{eqnarray}
(see reference \cite[(5.2) p.~86]{DoLa17} with Dondi). 
Then the following statement holds (see reference  \cite[Lem.~5.1, inequality at line 13 of p.~86]{DoLa17} with Dondi).
\begin{lem}
\label{lem:boest}
 Let ${\mathbf{a}}$ be as in (\ref{introd0}), (\ref{ellip}), (\ref{symr}).  Let $S_{ {\mathbf{a}} }$ be a fundamental solution of $P[{\mathbf{a}},D]$.  Let $\alpha\in]0,1]$. Let $\Omega$ be a bounded open subset of ${\mathbb{R}}^{n}$ of class $C^{1,\alpha}$. Then the following statements hold.
\begin{enumerate}
\item[(i)] If $\alpha\in]0,1[$, then
\begin{equation}
\label{lem:boest1}
b_{\Omega,\alpha}\equiv\sup
\biggl\{\biggr.
|x-y|^{n-1-\alpha}|\overline{B^{*}_{\Omega,y}}\left(S_{{\mathbf{a}}}(x-y)\right)|
:\, x,y\in \partial\Omega, x\neq y
\biggl.\biggr\}<+\infty\,.
\end{equation}
If $n>2$, then (\ref{lem:boest1}) holds also for $\alpha=1$. If $n=2$ and $DB_1(0)=0$, then (\ref{lem:boest1}) holds also for $\alpha=1$. 
\item[(ii)] If $n=2$ and $\alpha=1$, then
\begin{equation}
\label{boest1a}
b_{\Omega,\alpha}\equiv\sup
\biggl\{\biggr.
 \frac{|\overline{B^{*}_{\Omega,y}}\left(S_{{\mathbf{a}}}(x-y)\right)|	}{(1+|\ln |x-y||)}
:\, x,y\in \partial\Omega, x\neq y
\biggl.\biggr\}<+\infty.
\end{equation}
In particular, the kernel $\overline{B^{*}_{\Omega,y}}\left(S_{{\mathbf{a}}}(x-y)\right)$ belongs to
${\mathcal{K}}_{\epsilon,(\partial\Omega)\times (\partial\Omega)}$ for all $\epsilon\in]0,+\infty[$.
\item[(iii)]
\begin{eqnarray*}
\lefteqn{
\tilde{b}_{\Omega,\alpha}\equiv\sup
\biggl\{\biggr.\frac{|x'-y|^{n-\alpha}}{|x'-x''|}
|\overline{B^{*}_{\Omega,y}}\left(S_{{\mathbf{a}}}(x'-y)\right)
-
\overline{B^{*}_{\Omega,y}}\left(S_{{\mathbf{a}}}(x''-y)\right)|:\,
}
\\ \nonumber
&& \qquad\qquad
x',x''\in\partial\Omega, x'\neq x'', y\in\partial\Omega\setminus{\mathbb{B}}_{n}(x',2|x'-x''|)
\biggl.\biggr\}<+\infty\,.
\end{eqnarray*}
\end{enumerate}
\end{lem}
By applying equality (\ref{eq:boest3}),   we can  compute a formula for the  tangential gradient with respect to its first variable of the kernel of the double layer potential and establish some of its properties. 
To do so we introduce the following technical lemma 
(see reference  \cite[Lem. 3.2 (v), 3.3]{DoLa17} with Dondi).
\begin{lem}
\label{lem:fanes}
Let $Y$ be a nonempty bounded subset of ${\mathbb{R}}^{n}$. Then the following statements hold.  
\begin{enumerate}
\item[(i)] Let $F\in {\mathrm{Lip}}(\partial{\mathbb{B}}_{n}(0,1)\times [0,{\mathrm{diam}}\,(Y)])$ with
\begin{eqnarray*}
\lefteqn{
{\mathrm{Lip}}(F)
\equiv\biggl\{\biggr.
\frac{|F(\theta',r')-F(\theta'',r'')|}{ |\theta'-\theta''|+|r'-r''| }:\,
}
\\ \nonumber
&&\qquad 
(\theta',r'),(\theta'',r'')\in \partial{\mathbb{B}}_{n}(0,1)\times [0,{\mathrm{diam}}\,(Y)],\ (\theta',r')\neq (\theta'',r'')
\biggl.\biggr\}\,.
\end{eqnarray*}
Then 
\begin{eqnarray}
\label{lem:fanes1}
\lefteqn{
\left|
F\left(
\frac{x'-y}{|x'-y|},|x'-y|
\right)
-
F\left(
\frac{x''-y}{|x''-y|},|x''-y|
\right)
\right|
}
\\ \nonumber
&&\quad
\leq
{\mathrm{Lip}}(F) (2+   {\mathrm{diam}}\,(Y))
\frac{|x'-x''|}{|x'-y|}\,
\quad\forall y\in Y
\setminus {\mathbb{B}}_{n}(x',2|x'-x''|)\,,
\end{eqnarray}
for all $x',x''\in Y$, $x'\neq x''$. In particular, if $f\in C^{1}(\partial{\mathbb{B}}_{n}(0,1)\times{\mathbb{R}},{\mathbb{C}})$, then
\begin{eqnarray*}
\lefteqn{
M_{f,Y}\equiv\sup\biggl\{\biggr.
\left|
f\left(
\frac{x'-y}{|x'-y|},|x'-y|
\right)
-
f\left(
\frac{x''-y}{|x''-y|},|x''-y|
\right)
\right|\frac{|x'-y|}{|x'-x''|}
}
\\ \nonumber
&& \qquad\qquad\qquad\qquad  
:\,x',x''\in Y, x'\neq x'',  y\in Y   
\setminus {\mathbb{B}}_{n}(x',2|x'-x''|)
\biggl.\biggr\} 
\end{eqnarray*}
is finite and thus the kernel $f\left(
\frac{x-y}{|x-y|},|x-y|
\right)$ belongs to  ${\mathcal{K}}_{0,1,1}(Y\times Y)$.

\item[(ii)] Let $W$ be an open neighbourhood of $\overline{(Y-Y)}$. Let $f\in C^{1}(W,{\mathbb{C}})$. Then 
\begin{eqnarray*}
\lefteqn{
\tilde{M}_{f,Y}\equiv
 \sup\biggl\{\biggr.
|
f(x'-y)-f(x''-y)|\,|x'-x''|^{-1}:\,
}
\\ \nonumber
&&\qquad\qquad\qquad\qquad\qquad\qquad\quad
 x',x''\in Y, x'\neq x'',  y\in Y 
\biggl.\biggr\}<+\infty\,.
\end{eqnarray*}
Here $Y-Y\equiv\{y_{1}-y_{2}:\ y_{1}, y_{2}\in Y\}$. In particular, the kernel $f(x-y)$ belongs to the class ${\mathcal{K}}_{0,0,1}(Y\times Y)$, which is continuously imbedded into ${\mathcal{K}}_{0,1,1}(Y\times Y)$.
\item[(iii)] The kernel $\ln|x-y|$ belongs to ${\mathcal{K}}_{\epsilon,1,1}(Y\times Y)$ for all $\epsilon\in]0,1[$.
\end{enumerate}
\end{lem}
{\bf Proof.}  For the proof of (i), the first part of (ii) and (iii), we refer to the above mentioned paper  \cite[Lem. 3.2 (v), 3.3]{DoLa17}. The imbedding of the second part of (ii) follows by the imbedding Proposition \ref{prop:kerem} (ii). \hfill  $\Box$ 

\vspace{\baselineskip}
We are now ready to prove  the following statement. For the definition of tangential gradient $ {\mathrm{grad}}_{\partial\Omega}$ and tangential divergence $ {\mathrm{div}}_{\partial\Omega}$, we refer   to Kirsch and Hettlich \cite[A.5]{KiHe15}, Chavel~\cite[Chap.~1]{Cha84}. 
\begin{lem}\label{lem:tgdlgen}  
 Let ${\mathbf{a}}$ be as in (\ref{introd0}), (\ref{ellip}), (\ref{symr}).  Let $S_{ {\mathbf{a}} }$ be a fundamental solution of $P[{\mathbf{a}},D]$.   Let $\alpha\in]0,1]$. Let $\Omega$ be a bounded open subset of ${\mathbb{R}}^{n}$ of class $C^{1,\alpha}$.  Then the following statements hold.
 \begin{enumerate}
\item[(i)]  If  $h\in\{1,\dots,n\}$, then
\begin{eqnarray}\label{lem:tgdlgen1}
\lefteqn{
({\mathrm{grad}}_{\partial\Omega,x}\overline{B^{*}_{\Omega,y}}\left(S_{{\mathbf{a}}}(x-y)\right))_h
}
\\ \nonumber
&&\qquad
=\frac{\partial}{\partial x_h}\overline{B^{*}_{\Omega,y}}\left(S_{{\mathbf{a}}}(x-y)\right))
-
\nu_h(x)\sum_{l=1}^n \nu_l(x)
\frac{\partial}{\partial x_l}\overline{B^{*}_{\Omega,y}}\left(S_{{\mathbf{a}}}(x-y)\right)) 
\\ \nonumber
&& 
=\frac{n}{s_n\sqrt{\det a^{(2)} }}
 \frac{(x-y)^t\cdot\nu(y)}{|T^{-1}(x-y)|^n}
\\ \nonumber
&& \ 
\times
   \sum_{l=1}^n\nu_l(x)
 \biggl[
 \nu_l(x)
 \frac{
 \sum_{j,z=1}^n (T^{-1})_{jz}(x_z-y_z)(T^{-1})_{jh}
  }{|T^{-1}(x-y)|^2}  
  \\ \nonumber
&& \ 
 -
 \nu_h(x)\frac{\sum_{j,z=1}^n (T^{-1})_{jz}(x_z-y_z)(T^{-1})_{jl} 
  }{|T^{-1}(x-y)|^2}\biggr]
\\ \nonumber
&&\   
 -
 \frac{\sum_{l=1}^n\nu_l(x)
 \bigl[
 \nu_l(x) \nu_h(y)  
 -  \nu_h(x)  \nu_l(y) 
 \bigr]}{
 s_n\sqrt{\det a^{(2)} }   |T^{-1}(x-y)|^{n}} 
    \\ \nonumber
&&\   
-(2-n)|x-y|^{1-n}A_2\left(\frac{x-y}{|x-y|},|x-y|\right)a^{(2)}\nu(y)
\\ \nonumber
&&\   
\times
\sum_{l=1}^n\nu_l(x)\bigl[
\nu_l(x)\frac{x_h-y_h}{|x-y|}-\nu_h(x)\frac{x_l-y_l}{|x-y|}
\bigr]
\\ \nonumber
&&\   
-\sum_{j=1}^n\frac{\partial A_2}{\partial y_j}
\left(\frac{x-y}{|x-y|},|x-y|\right)a^{(2)}\nu(y)|x-y|^{-n}
\\ \nonumber
&&\   
\times
\sum_{l=1}^n\nu_l(x)
\biggl[
\nu_l(x)\biggl(
\delta_{jh}|x-y|-\frac{(x_j-y_j)(x_h-y_h)}{|x-y|}
\biggr)
\\ \nonumber
&&\   
-\nu_h(x)
\biggl(
\delta_{jl}|x-y|-\frac{(x_j-y_j)(x_l-y_l)}{|x-y|}
\biggr)
\biggr]
\\ \nonumber
&&\   
-\frac{\partial A_2}{\partial r}
\left(\frac{x-y}{|x-y|},|x-y|\right)a^{(2)}\nu(y)
\\ \nonumber
&&\    
\times
\sum_{l=1}^n\nu_l(x)
\left[
\nu_l(x)\frac{x_h-y_h}{|x-y|^{n-1}}-\nu_h(x)\frac{x_l-y_l}{|x-y|^{n-1}}
\right]
\\ \nonumber
&&\    
-\sum_{j,z=1}^n
\sum_{l=1}^n\nu_l(x)
\left[
\nu_l(x)
\frac{\partial^2B_1}{\partial x_h\partial x_j}(x-y)
-
\nu_h(x)
\frac{\partial^2B_1}{\partial x_l\partial x_j}(x-y)
\right]
\\ \nonumber
&&\    
\times
a_{jz}\nu_z(y) \ln |x-y|
\\ \nonumber
&&\   
-DB_1(x-y)a^{(2)}\nu(y) 
\sum_{l=1}^n\nu_l(x)
\left[
\nu_l(x)\frac{x_h-y_h}{|x-y|^2}-\nu_h(x)\frac{x_l-y_l}{|x-y|^2}
\right]
\\ \nonumber
&&\   
-\sum_{j,s	 =1}^n
\sum_{l=1}^n\nu_l(x)
\left[\nu_l(x)\frac{\partial^2C}{\partial x_h\partial x_j}(x-y)
-
\nu_h(x)\frac{\partial^2C}{\partial x_l\partial x_j}(x-y)
\right]
a_{ 
 js	 
}\nu_{   s	 	
}
(y) 
\\ \nonumber
&&\   
-\nu(y)^t\cdot a^{(1)} 
\sum_{l=1}^n\nu_l(x)
\left[\nu_l(x)\frac{\partial S_{{\mathbf{a}}} }{\partial x_h}(x-y)
- \nu_h(x)\frac{\partial S_{{\mathbf{a}}} }{\partial x_l}(x-y)
\right] 
 \end{eqnarray}
 for all $(x,y)\in (\partial\Omega)^2\setminus {\mathbb{D}}_{\partial\Omega}$, where we understand that the symbols
 \[
 \frac{\partial A_2}{\partial y_j}\qquad\forall j\in\{1,\dots,n\}
 \]
 denote    partial derivatives of any of the analytic extensions of $A_2$ to an open neighborhood of $(\partial{\mathbb{B}}_n(0,1))\times {\mathbb{R}}$ in ${\mathbb{R}}^{n+1}$.
 \item[(ii)] The kernel   $ {\mathrm{grad}}_{\partial\Omega,x}\overline{B^{*}_{\Omega,y}}\left(S_{{\mathbf{a}}}(x-y)\right) $
  belongs to $\left({\mathcal{K}}_{n-\alpha,n,\alpha}(\partial\Omega\times \partial\Omega)\right)^n$.
 \end{enumerate}
\end{lem}
{\bf Proof.}  (i) By formula (\ref{eq:boest3}), we have
\begin{eqnarray*}
\lefteqn{
\frac{\partial}{\partial x_h}\overline{B^{*}_{\Omega,y}}\left(S_{{\mathbf{a}}}(x-y)\right))
}
 \\ \nonumber
&&\qquad\ \
=-\frac{(-n)}{s_n\sqrt{\det a^{(2)} }}
 \sum_{j,z=1}^n\frac{(T^{-1})_{jz}(x_z-y_z)(T^{-1})_{jh} }{|T^{-1}(x-y)|^2} 
 \frac{(x-y)^t\cdot\nu(y)}{|T^{-1}(x-y)|^n}
 \\ \nonumber
&&\qquad\ \
 -\frac{1}{s_n\sqrt{\det a^{(2)} }}
 |T^{-1}(x-y)|^{-n}  \nu_h(y) 
 \\ \nonumber
&&\qquad\ \ 
 -(2-n)|x-y|^{1-n} \frac{x_h-y_h}{|x-y|} A_{2}(\frac{x-y}{|x-y|},|x-y|)a^{(2)}\nu(y)
\\ \nonumber
&& \qquad\ \ 
 -\sum_{j=1}^n\frac{\partial A_2}{\partial y_j}(\frac{x-y}{|x-y|},|x-y|)a^{(2)}\nu(y)
\frac{\delta_{jh}|x-y|-\frac{(x_j-y_j)(x_h-y_h)}{|x-y|}}{|x-y|^n}   
\\ \nonumber
&& \qquad\ \
- \frac{\partial A_2}{\partial r}(\frac{x-y}{|x-y|},|x-y|)a^{(2)}\nu(y)\frac{x_h-y_h}{|x-y|^{n-1}}  
\\ \nonumber
&& \qquad\ \ 
-\sum_{j,z=1}^n\frac{\partial^2B_1}{\partial x_h\partial x_j}(x-y)a_{jz}\nu_z(y) \ln |x-y| 
\\ \nonumber
&& \qquad\ \
-DB_1(x-y)a^{(2)}\nu(y) \frac{x_h-y_h}{|x-y|^2} 
\\ \nonumber
&& \qquad\ \ 
-\sum_{j,s=1}^n\frac{\partial^2C}{\partial x_h\partial x_j}(x-y)a_{js}\nu_s(y) 
-\nu(y)^t\cdot a^{(1)} \frac{\partial S_{{\mathbf{a}}} }{\partial x_h}(x-y) 
\end{eqnarray*} 
 for all $(x,y)\in (\partial\Omega)^2\setminus {\mathbb{D}}_{\partial\Omega}$. Then the definition of tangential gradient implies 
the validity of formula (\ref{lem:tgdlgen1}). 	
 
We now turn to the proof of (ii). If suffices to  show that if $ h\in\{1,\dots,n\} $, then
 each addendum in the right hand side of formula (\ref{lem:tgdlgen1}) belongs to the class ${\mathcal{K}}_{n-\alpha,n,\alpha}(\partial\Omega\times \partial\Omega)$.
 
 By Lemma \ref{lem:knnnone} the kernel $\frac{1}{ |T^{-1}(x-y)|^n}$ belongs to ${\mathcal{K}}_{n ,n+1,1}(\partial\Omega\times \partial\Omega)$.
 Since  there exists $c_{\Omega,\alpha}\in]0,+\infty[$ such that
\[
|\nu (y)\cdot (x-y)|\leq c_{\Omega,\alpha} |x-y|^{1+\alpha}
\qquad
\forall x,y\in \partial\Omega 
\]
 the kernel  $\nu(y)\cdot  (x-y)$ belongs to ${\mathcal{K}}_{-1-\alpha,-\alpha,1}(\partial\Omega\times \partial\Omega)$ (cf.~\textit{e.g.},  reference \cite[Lem.~3.4 and p.~87 line 8]{DoLa17} with Dondi). Then the product Theorem \ref{thm:kerpro} implies that the kernel 
$\frac{\nu(y)(x-y) }{|T^{-1}(x-y)|^{n}}$ belongs to ${\mathcal{K}}_{n-1-\alpha,n-\alpha,1}(\partial\Omega\times \partial\Omega)$.   By Lemma \ref{lem:kelem}, ${\mathcal{K}}_{n-1-\alpha,n-\alpha,1}(\partial\Omega\times \partial\Omega)$ is contained in ${\mathcal{K}}_{n-1-\alpha,n-1,\alpha}(\partial\Omega\times \partial\Omega)$.

By Lemma \ref{lem:knnnone} the kernel $\frac{x_h-y_h}{|T^{-1}(x-y)|^2}$ belongs to ${\mathcal{K}}_{1 ,2,1}(\partial\Omega\times \partial\Omega)$. By Lemma \ref{lem:kelem}, ${\mathcal{K}}_{1 ,2,1}(\partial\Omega\times \partial\Omega)$ is contained in ${\mathcal{K}}_{1 ,1+\alpha,\alpha}(\partial\Omega\times \partial\Omega)$. Then the   $\alpha$-H\"{o}lder continuity of $\nu$ and  Propostion \ref{prop:prkerho}  imply	 
 that 
\begin{eqnarray*}
\lefteqn{\sum_{l=1}^n\nu_l(x)
 \biggl[
 \nu_l(x)
 \frac{
 \sum_{j,z=1}^n (T^{-1})_{jz}(x_z-y_z)(T^{-1})_{jh}
  }{|T^{-1}(x-y)|^2}  
  }
  \\ \nonumber
&& \qquad \qquad \qquad \qquad 
 -
 \nu_h(x)\frac{\sum_{j,z=1}^n (T^{-1})_{jz}(x_z-y_z)(T^{-1})_{jl} 
  }{|T^{-1}(x-y)|^2}\biggr]
  \end{eqnarray*}
belongs to ${\mathcal{K}}_{1 ,1+\alpha,\alpha }(\partial\Omega\times \partial\Omega) $. Then the product Theorem \ref{thm:kerpro} (ii) implies that
 \begin{eqnarray}\label{lem:tgdlgen2}
 \lefteqn{
\frac{(x-y)^t\cdot\nu(y)}{|T^{-1}(x-y)|^n}
\sum_{l=1}^n\nu_l(x)
 \biggl[
 \nu_l(x)
 \frac{
 \sum_{j,z=1}^n (T^{-1})_{jz}(x_z-y_z)(T^{-1})_{jh}
  }{|T^{-1}(x-y)|^2}  
  }
  \\ \nonumber
&&\qquad\qquad\qquad\qquad\qquad\qquad
 -
 \nu_h(x)\frac{\sum_{j,z=1}^n (T^{-1})_{jz}(x_t-y_t)(T^{-1})_{jl} 
  }{|T^{-1}(x-y)|^2}\biggr]
  \\ \nonumber
&&\qquad
\in {\mathcal{K}}_{n-\alpha ,n,\alpha}(\partial\Omega\times \partial\Omega)
 \,.
\end{eqnarray}
We now consider the second addendum in the right hand side of formula (\ref{lem:tgdlgen1}) and we observe that
\begin{eqnarray*}
\lefteqn{
 \frac{\sum_{l=1}^n\nu_l(x)
 \bigl[
 \nu_l(x) \nu_h(y)  
 -  \nu_h(x)  \nu_l(y) 
 \bigr]}{
 s_n\sqrt{\det a^{(2)} }   |T^{-1}(x-y)|^{n}} 
}
\\ \nonumber
&&\qquad
=\frac{\sum_{l=1}^n\nu_l(x)
 \bigl[
 \nu_l(x) (\nu_h(y) -\nu_h(x))
 -  \nu_h(x) ( \nu_l(y) -\nu_l(x))
 \bigr]}{
 s_n\sqrt{\det a^{(2)} }   |T^{-1}(x-y)|^{n}} 
\end{eqnarray*} 
for all $(x,y)\in (\partial\Omega)^2\setminus {\mathbb{D}}_{\partial\Omega}$. 
   Since $\nu$ is $\alpha$-H\"{o}lder continuous, Lemma \ref{lem:hoker} implies that  $\nu_h(x)-\nu_h(y)$ belongs to ${\mathcal{K}}_{-\alpha,0,\alpha}(\partial\Omega\times \partial\Omega)$. By Lemma \ref{lem:knnnone} the kernel $\frac{1}{|T^{-1}(x-y)|^n}$ belongs to ${\mathcal{K}}_{n ,n+1,1}(\partial\Omega\times \partial\Omega)\subseteq 
{\mathcal{K}}_{n ,n+1-(1-\alpha),1-(1-\alpha)}(\partial\Omega\times \partial\Omega)$. Then the product Theorem \ref{thm:kerpro} (ii) implies that
 \[ 
\frac{\nu_h(x)-\nu_h(y)}{|T^{-1}(x-y)|^n}\in {\mathcal{K}}_{n-\alpha ,n+\alpha-\alpha,\alpha}(\partial\Omega\times \partial\Omega)
\,.
\]
Then the   $\alpha$-H\"{o}lder continuity of $\nu$ and  Propostion \ref{prop:prkerho} implies that
\[
\sum_{l=1}^n \frac{(\nu_l(x)-\nu_l(y))}{|T^{-1}(x-y)|^n}\nu_l(x)\nu_h(x)\in {\mathcal{K}}_{n-\alpha ,n  ,\alpha}(\partial\Omega\times \partial\Omega)\,.
\]
Hence,
\begin{equation}\label{lem:tgdlgen3}
\frac{\sum_{l=1}^n\nu_l(x)
 \bigl[
 \nu_l(x) \nu_h(y) 
 -  \nu_h(x) \nu_l(y) 
 \bigr]}{
    |T^{-1}(x-y)|^{n}} 
\in   {\mathcal{K}}_{n-\alpha ,n  ,\alpha}(\partial\Omega\times \partial\Omega)\,.
\end{equation}
We now consider the third addendum in the right hand side of formula (\ref{lem:tgdlgen1}). Since $A_2$ is real analytic in 
$\partial{\mathbb{B}}_n(0,1) \times {\mathbb{R}}$,  Lemma \ref{lem:fanes} (i) implies that
the kernel $A_2\left(\frac{x -y}{|x -y|},|x -y|\right)$ belongs to
${\mathcal{K}}_{0,1,1}(\partial\Omega\times\partial\Omega)$. Since the function 
$|\xi|^{1-n}\frac{\xi_h}{|\xi|}$ of the variable $\xi\in {\mathbb{R}}^n\setminus\{0\}$ is positively homogeneous of degree $-(n-1)$, Lemma \ref{lem:knnnone} implies that 
the kernel $|x-y|^{1-n}\frac{x_h-y_h}{|x-y|}$ is of class ${\mathcal{K}}_{n-1,n,1}(\partial\Omega\times\partial\Omega)$. Then the product Theorem \ref{thm:kerpro} (ii) and Proposition \ref{prop:prkerho} (iii)
imply that the kernel 
\[
-(2-n)|x-y|^{1-n} \frac{x_h-y_h}{|x-y|} A_2\left(\frac{x-y}{|x-y|},|x-y|\right)a^{(2)}\nu(y)
\]
belongs to the class ${\mathcal{K}}_{n-1,n,1}(\partial\Omega\times\partial\Omega)$. By the imbedding Proposition \ref{prop:kerem} (ii) with
\[
s_1=n-1\,,\quad s_2=n\,,\quad s_3=1\,,\qquad
t_1=n-\alpha\,,\quad t_2=n\,,\quad t_3=\alpha\,,
\]
 ${\mathcal{K}}_{n-1,n,1}(\partial\Omega\times\partial\Omega)$
is contained in ${\mathcal{K}}_{n-\alpha,n,\alpha}(\partial\Omega\times\partial\Omega)$. 
Since the components of $\nu$ are of class $C^{0,\alpha}$, the product Proposition \ref{prop:prkerho}  (ii) implies that
\begin{eqnarray}\label{lem:tgdlgen4}
\lefteqn{-(2-n)|x-y|^{1-n}A_2\left(\frac{x-y}{|x-y|},|x-y|\right)a^{(2)}\nu(y)
}
\\ \nonumber
&&\qquad
\sum_{l=1}^n\nu_l(x)\bigl[
\nu_l(x)\frac{x_h-y_h}{|x-y|}-\nu_h(x)\frac{x_l-y_l}{|x-y|}
\bigr]
\in {\mathcal{K}}_{n-\alpha,n,\alpha}(\partial\Omega\times\partial\Omega)\,.
\end{eqnarray}
 We now consider the fourth addendum in the right hand side of formula (\ref{lem:tgdlgen1}). Let $j\in \{1,\dots,n\}$. Since $\frac{\partial A_2}{\partial y_j}$ is real analytic in 
$\partial{\mathbb{B}}_n(0,1) \times {\mathbb{R}}$,   Lemma \ref{lem:fanes} (i) implies that
the kernel $\frac{\partial A_2}{\partial y_j}\left(\frac{x -y}{|x -y|},|x -y|\right)$ belongs to
${\mathcal{K}}_{0,1,1}(\partial\Omega\times\partial\Omega)$. By Lemma \ref{lem:kelem}, 
\[
{\mathcal{K}}_{0,1,1}(\partial\Omega\times\partial\Omega)\subseteq
{\mathcal{K}}_{0,1-(1-\alpha),1-(1-\alpha)}(\partial\Omega\times\partial\Omega)
={\mathcal{K}}_{0, \alpha , \alpha}(\partial\Omega\times\partial\Omega)
\,.
\]
Since the functions 
$|\xi|^{-(n-1)}$ and $|\xi|^{-n-1}\xi_j \xi_l$  of the variable $\xi\in {\mathbb{R}}^n\setminus\{0\}$ are positively homogeneous of degree $-(n-1)$, Lemma \ref{lem:knnnone} implies that 
the kernels $|x-y|^{-(n-1)}$ and $|x-y|^{-n-1}(x_j-y_j)(x_l-y_l)$ are of class ${\mathcal{K}}_{n-1,n,1}(\partial\Omega\times\partial\Omega)$. By Lemma \ref{lem:kelem}, ${\mathcal{K}}_{n-1,n,1}(\partial\Omega\times\partial\Omega)$ is contained in ${\mathcal{K}}_{n-1,n-1+\alpha,\alpha}(\partial\Omega\times\partial\Omega)$.
Then the product Theorem \ref{thm:kerpro} (ii) implies that the product is continuous from
\[
  {\mathcal{K}}_{n-1,n-1+\alpha,\alpha}(\partial\Omega\times\partial\Omega)
\times {\mathcal{K}}_{0, \alpha , \alpha}(\partial\Omega\times\partial\Omega)
\quad\text{to}\quad {\mathcal{K}}_{n-1,n-1+\alpha,\alpha}(\partial\Omega\times\partial\Omega)\,.
\]
Then the $\alpha$-H\"{o}lder continuity of the components of $\nu$, Proposition \ref{prop:prkerho} (ii), (iii) and the imbedding Proposition \ref{prop:kerem} (iii) imply that  
\begin{eqnarray}\label{lem:tgdlgen5}
\lefteqn{
-\sum_{j=1}^n\frac{\partial A_2}{\partial y_j}
\left(\frac{x-y}{|x-y|},|x-y|\right)a^{(2)}\nu(y)|x-y|^{-n}
}
\\ \nonumber
&&\quad 
\times
\sum_{l=1}^n\nu_l(x)
\biggl[
\nu_l(x)\biggl(
\delta_{jh}|x-y|-\frac{(x_j-y_j)(x_h-y_h)}{|x-y|}
\biggr)
\\ \nonumber
&&\quad 
-\nu_h(x)
\biggl(
\delta_{jl}|x-y|-\frac{(x_j-y_j)(x_l-y_l)}{|x-y|}
\biggr)
\biggr]
\\ \nonumber
&&\quad 
\in {\mathcal{K}}_{n-1,n-1+\alpha,\alpha}(\partial\Omega\times\partial\Omega)
\subseteq {\mathcal{K}}_{n-\alpha,n,\alpha}(\partial\Omega\times\partial\Omega)\,.
\end{eqnarray}
 We now consider the fifth addendum in the right hand side of formula (\ref{lem:tgdlgen1}). Since $\frac{\partial A_2}{\partial r}$ is real analytic in 
$\partial{\mathbb{B}}_n(0,1) \times {\mathbb{R}}$,   Lemma \ref{lem:fanes} (i) implies that
the kernel $\frac{\partial A_2}{\partial r}\left(\frac{x -y}{|x -y|},|x -y|\right)$ belongs to
${\mathcal{K}}_{0,1,1}(\partial\Omega\times\partial\Omega)$ that is contained in $  {\mathcal{K}}_{0, \alpha, \alpha}(\partial\Omega\times\partial\Omega)$ (cf.~Lemma \ref{lem:kelem}). Since the function    $|\xi|^{-(n-1)} \xi_l$  of the variable $\xi\in {\mathbb{R}}^n\setminus\{0\}$ is positively homogeneous of degree $n-2$, Lemma \ref{lem:knnnone} implies that 
the kernels $|x-y|^{-(n-1)}(x_l-y_l)$  are of class ${\mathcal{K}}_{n-2,n-1,1}(\partial\Omega\times\partial\Omega)$, that is contained in ${\mathcal{K}}_{n-2,n-2+\alpha,\alpha}(\partial\Omega\times\partial\Omega)$ (cf.~Lemma \ref{lem:kelem}).
Then the product Theorem \ref{thm:kerpro} (ii)  implies that the product is continuous from
\[
  {\mathcal{K}}_{n-2,n-2+\alpha,\alpha}(\partial\Omega\times\partial\Omega)
\times {\mathcal{K}}_{0, \alpha , \alpha}(\partial\Omega\times\partial\Omega)
\quad\text{to}\quad {\mathcal{K}}_{n-2,n-2+\alpha,\alpha}(\partial\Omega\times\partial\Omega)\,.
\]
Then the  $\alpha$-H\"{o}lder continuity of the components of $\nu$, Proposition \ref{prop:prkerho} (ii), (iii) and the imbedding Proposition  \ref{prop:kerem} (iii) imply that  
\begin{eqnarray}\label{lem:tgdlgen6}
\lefteqn{
\frac{\partial A_2}{\partial r}
\left(\frac{x-y}{|x-y|},|x-y|\right)a^{(2)}\nu(y)
}
\\ \nonumber
&&\quad 
\times
\sum_{l=1}^n\nu_l(x)
\left[
\nu_l(x)\frac{x_h-y_h}{|x-y|^{n-1}}-\nu_h(x)\frac{x_l-y_l}{|x-y|^{n-1}}
\right]
\\ \nonumber
&&\quad
\in 
{\mathcal{K}}_{n-2,n-2+\alpha,\alpha}(\partial\Omega\times\partial\Omega)
\subseteq
{\mathcal{K}}_{n-\alpha,n,\alpha}(\partial\Omega\times\partial\Omega)\,.
\end{eqnarray}
We now consider the sixth addendum in the right hand side of formula (\ref{lem:tgdlgen1}). Since $B_1$ is analytic, Lemma \ref{lem:fanes} (ii) implies that the kernel $\frac{\partial^2B_1}{\partial x_l\partial x_j}(x-y)$ belongs to
${\mathcal{K}}_{0,1,1}(\partial\Omega\times\partial\Omega)$ that is contained in $  {\mathcal{K}}_{0, \alpha, \alpha}(\partial\Omega\times\partial\Omega)$ for each $j,l\in\{1,\dots,n\}$ (cf.~Lemma \ref{lem:kelem}).
Then the  $\alpha$-H\"{o}lder continuity of the components of $\nu$ and the product Proposition \ref{prop:prkerho} (ii), (iii)  imply that
\begin{eqnarray*}
\lefteqn{
\sum_{j,t=1}^n
\sum_{l=1}^n\nu_l(x)
\left[
\nu_l(x)
\frac{\partial^2B_1}{\partial x_h\partial x_j}(x-y)
-
\nu_h(x)
\frac{\partial^2B_1}{\partial x_l\partial x_j}(x-y)
\right]
a_{jt}\nu_t(y)
}
\\ \nonumber
&&\qquad\qquad\qquad\qquad\qquad\qquad\qquad\qquad\qquad\qquad
\in  {\mathcal{K}}_{0,\alpha,\alpha}(\partial\Omega\times\partial\Omega)
\,.
\end{eqnarray*}
 By Lemma \ref{lem:fanes} (iii) and by Lemma \ref{lem:kelem}, we have 
 \[
 \ln|x-y|\in {\mathcal{K}}_{\epsilon,1,1}(\partial\Omega\times\partial\Omega)
 \subseteq {\mathcal{K}}_{\epsilon,\alpha,\alpha}(\partial\Omega\times\partial\Omega)
 \qquad\forall\epsilon\in]0,1[\,.
 \]
Theorem \ref{thm:kerpro} (ii) implies that the product is continuous from
\[
{\mathcal{K}}_{0,\alpha,\alpha}(\partial\Omega\times\partial\Omega)
\times
{\mathcal{K}}_{\epsilon,\alpha,\alpha}(\partial\Omega\times\partial\Omega)
\quad\text{to}\quad
{\mathcal{K}}_{\epsilon,\alpha+\epsilon,\alpha}(\partial\Omega\times\partial\Omega)\,.
\]
Hence, inequalities $n-\alpha\geq \epsilon$, $\alpha\leq \alpha$ 
 and the imbedding Proposition  \ref{prop:kerem} (iii) 
 imply that
 \begin{eqnarray}\label{lem:tgdlgen7}
\lefteqn{
\sum_{j,z=1}^n
\sum_{l=1}^n\nu_l(x)
\left[
\nu_l(x)
\frac{\partial^2B_1}{\partial x_h\partial x_j}(x-y)
-
\nu_h(x)
\frac{\partial^2B_1}{\partial x_l\partial x_j}(x-y)
\right]	}
\\ \nonumber
&&\quad 
\times
a_{jz}\nu_z(y) \ln |x-y|
\in{\mathcal{K}}_{\epsilon,\alpha+\epsilon,\alpha}(\partial\Omega\times\partial\Omega)
\subseteq
{\mathcal{K}}_{n-\alpha,n,\alpha}(\partial\Omega\times\partial\Omega)\,.
\end{eqnarray}
 We now consider the seventh addendum in the right hand side of formula (\ref{lem:tgdlgen1}). Since $B_1$ is analytic, Lemma \ref{lem:fanes} (ii) and the product Proposition \ref{prop:prkerho} (iii) imply that
 $DB_1(x-y)a^{(2)}\nu(y) $ belongs to ${\mathcal{K}}_{0,1,1}(\partial\Omega\times\partial\Omega)$
 that is contained in $  {\mathcal{K}}_{0, \alpha, \alpha}(\partial\Omega\times\partial\Omega)$ (cf.~Lemma \ref{lem:kelem}).
 Since the functions    $|\xi|^{-2} \xi_l$  of the variable 
 $\xi\in {\mathbb{R}}^n\setminus\{0\}$ are positively homogeneous of degree $-1$, Lemma \ref{lem:knnnone} implies that 
the kernels $|x-y|^{-2}(x_l-y_l)$  are of class ${\mathcal{K}}_{1,2,1}(\partial\Omega\times\partial\Omega)$  that is contained in $  {\mathcal{K}}_{1,1+ \alpha, \alpha}(\partial\Omega\times\partial\Omega)$ (cf.~Lemma \ref{lem:kelem}).
 Hence the $\alpha$-H\"{o}lder continuity of the components of $\nu$ and the product Proposition \ref{prop:prkerho} (ii) imply that
\[
\sum_{l=1}^n\nu_l(x)
\left[
\nu_l(x)\frac{x_h-y_h}{|x-y|^2}-\nu_h(x)\frac{x_l-y_l}{|x-y|^2}
\right]\in {\mathcal{K}}_{1,1+ \alpha, \alpha}(\partial\Omega\times\partial\Omega).
\]
Theorem \ref{thm:kerpro} (ii) implies that the product is continuous from
\[
{\mathcal{K}}_{0,\alpha,\alpha}(\partial\Omega\times\partial\Omega)
\times
{\mathcal{K}}_{1,1+ \alpha, \alpha}(\partial\Omega\times\partial\Omega)
\quad\text{to}\quad
{\mathcal{K}}_{1,1+\alpha,\alpha}(\partial\Omega\times\partial\Omega) 
\]
and thus the  imbedding Proposition \ref{prop:kerem} (iii) implies that
 \begin{eqnarray}\label{lem:tgdlgen8}
\lefteqn{
DB_1(x-y)a^{(2)}\nu(y) 
\sum_{l=1}^n\nu_l(x)
\left[
\nu_l(x)\frac{x_h-y_h}{|x-y|^2}-\nu_h(x)\frac{x_l-y_l}{|x-y|^2}
\right]
}
\\ \nonumber
&&\qquad\qquad\qquad\qquad
\in {\mathcal{K}}_{1,1+\alpha,\alpha}(\partial\Omega\times\partial\Omega)
\subseteq
{\mathcal{K}}_{n-\alpha,n,\alpha}(\partial\Omega\times\partial\Omega)\,.
\end{eqnarray}
We now consider the eighth addendum in the right hand side of formula (\ref{lem:tgdlgen1}). Since $C$ is analytic, Lemma \ref{lem:fanes} (ii) implies that the kernel
$\frac{\partial^2C}{\partial x_l\partial x_j}(x-y)$ belongs to
${\mathcal{K}}_{0,1,1}(\partial\Omega\times\partial\Omega)$ that is contained in $  {\mathcal{K}}_{0, \alpha, \alpha}(\partial\Omega\times\partial\Omega)$ for each $j,l\in\{1,\dots,n\}$ (cf.~Lemma \ref{lem:kelem}). Then the  $\alpha$-H\"{o}lder continuity of the components of $\nu$, the product Proposition \ref{prop:prkerho} (ii), (iii) and the imbedding Proposition  \ref{prop:kerem} (iii) imply that
\begin{eqnarray}\label{lem:tgdlgen9}
\lefteqn{
\sum_{
 j,s	 
=1}^n
\sum_{l=1}^n\nu_l(x)
\left[\nu_l(x)\frac{\partial^2C}{\partial x_h\partial x_j}(x-y)
-
\nu_h(x)\frac{\partial^2C}{\partial x_l\partial x_j}(x-y)
\right]
a_{j
  s	 
}\nu_{
 s	 	
}(y) 
}
\\ \nonumber
&&\qquad\qquad\qquad\qquad\qquad\qquad 
\in  {\mathcal{K}}_{0,\alpha,\alpha}(\partial\Omega\times\partial\Omega)
\subseteq
{\mathcal{K}}_{n-\alpha,n,\alpha}(\partial\Omega\times\partial\Omega)\,.
\end{eqnarray}
We now consider the nineth addendum in the right hand side of formula (\ref{lem:tgdlgen1}). By reference  \cite[Rmk.~6.1]{DoLa17} with Dondi
the kernels $\frac{\partial S_{{\mathbf{a}}} }{\partial x_l}(x-y)$ belong to the class ${\mathcal{K}}_{n-1,n,1}(\partial\Omega\times\partial\Omega)$ that is contained in $  {\mathcal{K}}_{n-1,n-1 + \alpha, \alpha}(\partial\Omega\times\partial\Omega)$ for each $l\in\{1,\dots,n\}$ (cf.~Lemma \ref{lem:kelem}). 
Hence the $\alpha$-H\"{o}lder continuity of the components of $\nu$,  the product Proposition \ref{prop:prkerho} (ii), (iii) and the imbedding Proposition  \ref{prop:kerem} (iii) imply that
\begin{eqnarray}\label{lem:tgdlgen10}
\lefteqn{
-\nu(y)^t\cdot a^{(1)} 
\sum_{l=1}^n\nu_l(x)
\left[\nu_l(x)\frac{\partial S_{{\mathbf{a}}} }{\partial x_h}(x-y)
- \nu_h(x)\frac{\partial S_{{\mathbf{a}}} }{\partial x_l}(x-y)
\right] 
}
\\ \nonumber
&&\qquad\qquad\qquad\qquad
\in {\mathcal{K}}_{n-1,n-1 + \alpha, \alpha}(\partial\Omega\times\partial\Omega)
\subseteq
{\mathcal{K}}_{n-\alpha,n,\alpha}(\partial\Omega\times\partial\Omega)\,.
\end{eqnarray}
By   the memberships of (\ref{lem:tgdlgen2})--(\ref{lem:tgdlgen10}), we conclude that
 each addendum in the right hand side of formula (\ref{lem:tgdlgen1}) belongs to the class ${\mathcal{K}}_{n-\alpha,n,\alpha}(\partial\Omega\times \partial\Omega)$ and thus the proof is complete.\hfill  $\Box$ 

\vspace{\baselineskip}

\section{Continuity properties of the double layer potential}
As a consequence of Lemmas \ref{lem:boest} and \ref{lem:tgdlgen}, we can apply Theorem \ref{thm:iokreg} and prove the following classical result on the continuity of the double layer potential on the boundary (see Miranda \cite[15.VI]{Mi70}, where the author mentions a result of Giraud \cite{Gi32}. For the Laplace operator in case $n=2$ see Fichera and De Vito \cite[LXXXIII]{FiDe70}).
\begin{thm}\label{thm:dllregen}
 Let $n\in {\mathbb{N}}\setminus\{0,1\}$. 
 Let ${\mathbf{a}}$ be as in (\ref{introd0}), (\ref{ellip}), (\ref{symr}).  Let $S_{ {\mathbf{a}} }$ be a fundamental solution of $P[{\mathbf{a}},D]$.   Let $\alpha\in]0,1[$, $\beta\in]0,1]$, $\alpha+\beta>1$. 
 
  
 Let $\Omega$ be a bounded open subset of ${\mathbb{R}}^{n}$ of class $C^{1,\alpha}$.   Then the following statements hold.
 \begin{enumerate}
\item[(i)] If $\beta<1$, then the    operator $W_\Omega[{\mathbf{a}},S_{{\mathbf{a}}}   ,\cdot]$ from 
$C^{0,\beta}(\partial\Omega)$ to $C^{1,\alpha+\beta-1}(\partial\Omega)$ defined by (\ref{introd3})
 for all $\mu\in C^{0,\beta}(\partial\Omega)$ is linear and continuous.
\item[(ii)] If $\beta=1$, then the  operator $W_\Omega[{\mathbf{a}},S_{{\mathbf{a}}}   ,\cdot]$
 from 
$C^{0,\beta}(\partial\Omega)=C^{0,1}(\partial\Omega)$ to $C^{1,\omega_{\alpha+\beta-1}}(\partial\Omega)=C^{1,\omega_\alpha}(\partial\Omega)$ defined by (\ref{introd3})
 for all $\mu\in  C^{0,1}(\partial\Omega)$ 
 is linear and continuous.
\end{enumerate}
\end{thm} 
{\bf Proof.}  By formula (\ref{eq:boest3}), we have $
\overline{B^{*}_{\Omega,y}}\left(S_{{\mathbf{a}}}(\cdot-y)\right))\in C^1((\partial\Omega)\setminus\{y\})$ for all $y\in \partial\Omega$.	 
By Lemmas \ref{lem:boest} and \ref{lem:tgdlgen}, we know that the kernel of the double layer potential belongs to   ${\mathcal{K}}_{n-1-\alpha, n-\alpha,1}(\partial\Omega\times\partial\Omega)$ and that its tangential gradient with respect to the variable $x$ belongs to   $({\mathcal{K}}_{n-\alpha,n,\alpha}(\partial\Omega\times\partial\Omega))^n$. We now plan to apply Theorem \ref{thm:iokreg} (iii). We first note that 
reference \cite[Thm 9.2]{DoLa17} with Dondi implies that 
$W_\Omega[{\mathbf{a}},S_{{\mathbf{a}}}   ,1]\in C^{1,\alpha}(\partial\Omega)$. Moreover, \begin{eqnarray*}
&&\beta\leq 1\leq n-1<n-\alpha\equiv t_1=(n-1)+(1-\alpha)<(n-1)+\beta\,,
\\
&&t_2\equiv n\geq (n-1)+\beta\,,
\quad
0\leq s_1\equiv (n-1)-\alpha<n-1\,.
\end{eqnarray*}
(i) If $\beta<1$, then $t_2-\beta=n-\beta=(n-1)+1-\beta>n-1$, 
\[
\beta\leq 2\leq  t_2=n< n+\alpha+\beta-1=(n-1)+\beta+t_3\,,\ \ \text{where}\ t_3\equiv \alpha\,,
\]
and
\begin{eqnarray*}
\lefteqn{
\min\{\beta,(n-1)+\beta-t_1,(n-1)+t_3+\beta-t_2\}
}
\\ \nonumber
&& 
=\min\{\beta,(n-1)+\beta-(n-\alpha),(n-1)+\alpha+\beta-n\}
=\alpha+\beta-1\leq\alpha\,.
\end{eqnarray*}
Then
\begin{eqnarray*}
\lefteqn{
W_\Omega[{\mathbf{a}},S_{{\mathbf{a}}}   ,1]\in C^{1,\alpha}(\partial\Omega)
}
\\ \nonumber
&&\qquad
\subseteq C^{1,\alpha+\beta-1} (\partial\Omega)
=C^{1,\min\{\beta,(n-1)+\beta-t_1,(n-1)+t_3+\beta-t_2\}} (\partial\Omega)
\end{eqnarray*}
and Theorem \ref{thm:iokreg} (iii) (c) implies that $W_\Omega[{\mathbf{a}},S_{{\mathbf{a}}}   ,\cdot]$ is linear and continuous from $C^{0,\beta}(\partial\Omega)$ to
\[
C^{1,\min\{\beta,(n-1)+\beta-t_1,(n-1)+t_3+\beta-t_2\}}(\partial\Omega)
=C^{1,\alpha+\beta-1}(\partial\Omega)\,.
\]

(ii) If $\beta=1$, then $t_2-\beta=n-\beta=n-1$ and   
\[
C^{1,\max\{r^\beta, r^{(n-1)+\beta-t_1},\omega_{t_3}(\cdot)\}
}(\partial\Omega)
=C^{1,\max\{r,r^{\alpha},\omega_{\alpha}(\cdot)\}
} (\partial\Omega)
=C^{1,\omega_{\alpha}(\cdot)
} (\partial\Omega)\,.
\]
Then
\[
W_\Omega[{\mathbf{a}},S_{{\mathbf{a}}}   ,1]\in C^{1,\alpha}(\partial\Omega)\subseteq C^{1, \omega_{\alpha}(\cdot) }(\partial\Omega) 
=C^{1,
\max\{r^\beta,r^{(n-1)+\beta-t_1},\omega_{t_3}(\cdot)\}
}(\partial\Omega)
\]
and Theorem \ref{thm:iokreg} (iii) (cc) implies that $W_\Omega[{\mathbf{a}},S_{{\mathbf{a}}}   ,\cdot]$ is linear and continuous from $C^{0,\beta}(\partial\Omega)=C^{0,1}(\partial\Omega)$ to
\[
C^{1,
\max\{r^\beta,r^{(n-1)+\beta-t_1},\omega_{t_3}(\cdot)\}
}(\partial\Omega)
=C^{1, \omega_{\alpha}(\cdot) }(\partial\Omega) 
\]
 and thus the proof is complete.\hfill  $\Box$ 

\vspace{\baselineskip}   

 Next we introduce the following two technical  statements  in case $n=2$.
\begin{lem}\label{lem:c5o}
 Let $\Omega$ be a bounded open Lipschitz subset of ${\mathbb{R}}^{2}$. Then
 \[
 c_\Omega^{(v)}\equiv\sup_{x\in\partial\Omega,s\in]0,1/e[}|s\log s|^{-1}
 \int_{(\partial\Omega)\cap {\mathbb{B}}_2(0,s)}|\log|x-y||\,d\sigma_y<+\infty\,.
 \]
\end{lem}{\bf Proof.} By the Lemma  of the uniform cylinders, there exist  $r$, $\delta\in]0,1/e[$ such that if  $x\in \partial\Omega$, then there exist a $2\times 2$ orthogonal matrix
$R_x$ such that 
\[
C(x,R_x,r,\delta)\equiv x+ R_x^{t}({\mathbb{B}}_{2-1}(0,r)\times]-\delta,\delta[    ) 
\]
 is a coordinate cylinder for $ \Omega$ around $x$, \textit{i.e.}, there exists $\gamma_x\in C^{0,1}(\overline{{\mathbb{B}}_{1}(0,r)})$ such that \begin{eqnarray}
\label{prelim.cocylind1}
\lefteqn{R_x( \Omega -x )\cap ({\mathbb{B}}_{2-1}(0,r)\times ]-\delta,\delta[)
}
\\   \nonumber
&&\qquad\ 
=\left\{
(\eta,y)\in 
{\mathbb{B}}_{2-1}(0,r)\times ]-\delta,\delta[:\, y<\gamma_x(\eta)
\right\}
\equiv{\mathrm{hypograph}}_{s}(\gamma_x) 
\,,
\\   \nonumber
&&|\gamma_x(\eta)|<\delta/2\qquad\forall \eta\in {\mathbb{B}}_{2-1}(0,r)\,,\qquad
\gamma_x(0)=0\,,
\end{eqnarray}
  and the corresponding function $\gamma_x$ satisfies  the inequality
\[
 A\equiv\sup_{x\in \partial\Omega}\|\gamma_x\|_{ C^{0,1}(\overline{{\mathbb{B}}_{1}(0,r)}) }<+\infty 
 \]
 (cf.~\cite[Defn.~10.1, Lem.~10.1]{La20}).  By the continuity of the logarithm, it suffices to show that the supremum of the statement is finite with $s\in]0,r[$ and we note that $(\partial\Omega)\cap {\mathbb{B}}_2(x,s)\subseteq 
(\partial\Omega)\cap C(x,R_x,r,\delta)$ for all $s\in]0,r[$. Then we have
\begin{eqnarray*}
\lefteqn{
\int_{(\partial\Omega)\cap {\mathbb{B}}_2(x,s)}|\log|x-y||\,d\sigma_y
}
\\ \nonumber
&&\quad
\leq \int_{\{\eta\in ]-r,r[:|\eta|^2+\gamma_x(\eta)^2< s^2\}}
|\log|(\eta,\gamma_x(\eta))| \,|\,d\eta\sqrt{1+
 {\mathrm{ess}}	 
\sup|  \gamma'_x	 
|^2}
\\ \nonumber
&&\quad
\leq \int_{\{\eta\in ]-r,r[:|\eta|< s\}}
|\log|\eta| \,|\,d\eta\sqrt{1+A^2}
\leq 2\left[\eta-\eta\log\eta\right]_{\eta=0^+}^{\eta=s}\sqrt{1+A^2}
\\ \nonumber
&&\quad
\leq 4 |s\log s|\sqrt{1+A^2}\qquad\forall x\in\partial\Omega,s\in]0,1/e[\,.
\end{eqnarray*} \hfill  $\Box$ 

\vspace{\baselineskip}

\begin{prop}
\label{prop:v0a} 
Let $n=2$. Let ${\mathbf{a}}$ be as in (\ref{introd0}), (\ref{ellip}). Let $S_{ {\mathbf{a}} }$ be a fundamental solution of $P[{\mathbf{a}},D]$. 
  Let $\Omega$ be a bounded open  Lipschitz subset of ${\mathbb{R}}^{2}$. Let $S_{ {\mathbf{a}} }$ be a fundamental solution of $P[{\mathbf{a}},D]$.  Then   $v_\Omega[ S_{ {\mathbf{a}} }   ,\cdot   ]$ is continuous from $L^{\infty}(\partial\Omega)$ to $C^{0,\omega_{1}(\cdot)  }(\partial\Omega)$.
\end{prop}
{\bf Proof.} By Theorem 7.2 of \cite{DoLa17} we already know that $v_\Omega[ S_{ {\mathbf{a}} }   ,\cdot   ]$ is continuous from $L^{\infty}(\partial\Omega)$ to $C^{0 }(\partial\Omega)$.  We now take  $\mu\in L^\infty(\partial\Omega)$ and  we turn to estimate the H\"{o}lder constant  of $v_\Omega[ S_{ {\mathbf{a}} }   ,\mu  ]$. By formula (\ref{prop:ourfs1}) above, 
 by the inequality $|T^{-1}x|\geq |T|^{-1}|x|$ for $x\in{\mathbb{R}}^2\setminus\{0\}$ and by Lemma 4.2 (ii) of \cite{DoLa17}, there exists a constant $c\in]0,+\infty[$ such that
\begin{eqnarray*}
\lefteqn{  |\log|\xi||^{-1}|S_{ {\mathbf{a}} }(\xi)|\leq c\quad\forall   \xi\in
{\mathbb{B}}_2(0,1/e)\setminus\{0\}
\,,
}
\\
\lefteqn{ \frac{|x'-y| }{|x'-x''|}
|S_{ {\mathbf{a}} }(x'-y)-
S_{ {\mathbf{a}} }(x''-y)
|
\leq c
}
\\
&&\qquad\qquad
\forall x',x''\in \partial\Omega, x'\neq x'',  y\in (\partial\Omega)
\setminus {\mathbb{B}}_{n}(x',2|x'-x''|)\,.
\end{eqnarray*}
 Let $x',x''\in \partial\Omega$, $x'\neq x''$.   By Remark \ref{rem:om4}, there is  no loss of generality in assuming that $0<3|x'-x''|\leq 1/e$. Then the inclusion
${\mathbb{B}}_{2}(x',2|x'-x''|)\subseteq {\mathbb{B}}_{2}(x'',3|x'-x''|)$ and the triangular inequality imply  that
\begin{eqnarray}
\label{prop:k0a2}
\lefteqn{
|v_\Omega[ S_{ {\mathbf{a}} }   ,\mu  ](x')-v_\Omega[ S_{ {\mathbf{a}} }   ,\mu  ](x'')|
}
\\ \nonumber
&&\qquad
\leq   \|\mu\|_{  L^{\infty}(\partial\Omega)  }
\biggl\{\biggr.
\int_{{\mathbb{B}}_{2}(x',2|x'-x''|)\cap\partial\Omega}|S_{ {\mathbf{a}} }(x'-y)|\,d\sigma_{y}
 \\
\nonumber
&&\qquad\quad
+
 \int_{{\mathbb{B}}_{2}(x'',3|x'-x''|)\cap\partial\Omega}|S_{ {\mathbf{a}} }(x''-y)|\,d\sigma_{y}
 \\
\nonumber
&&\qquad\quad
+ 
\int_{\partial\Omega\setminus {\mathbb{B}}_{2}(x',2|x'-x''|)}
 |\, S_{ {\mathbf{a}} }(x'-y)-S_{ {\mathbf{a}} }(x''-y)  \,|
 \,d\sigma_{y}\biggl.\biggr\} \,.
\end{eqnarray} 
Then Lemma \ref{lem:c5o}  implies that
\begin{eqnarray}
\label{prop:k0a3}
\lefteqn{
\int_{{\mathbb{B}}_{2}(x',2|x'-x''|)\cap\partial\Omega}|S_{ {\mathbf{a}} } (x'-y)|\,d\sigma_{y}
}
\\ \nonumber
&& \qquad\qquad\qquad\qquad\qquad\qquad 
 +
 \int_{{\mathbb{B}}_{2}(x'',3|x'-x''|)\cap\partial\Omega}|S_{ {\mathbf{a}} } (x''-y)|\,d\sigma_{y}
\\ \nonumber
&& \qquad 
\leq
c
\biggl\{\biggr.
\int_{    {\mathbb{B}}_{2}(x',2|x'-x''|)\cap\partial\Omega     }
  | \log|x'-y|| d\sigma_{y}       
\\ \nonumber
&& \qquad\qquad\qquad\quad\quad
+
 \int_{   {\mathbb{B}}_{2}(x'',3|x'-x''|)\cap\partial\Omega    }
|\log|x''-y|| d\sigma_{y}      
 \biggl.\biggr\}
\\ \nonumber
&&  \qquad\qquad\qquad\quad
\leq
c2c^{(v)}_{\Omega}3|x'-x''||\log(3|x'-x''|)|
\\ \nonumber
&&  \qquad\qquad\qquad\quad
\leq 6cc^{(v)}_{\Omega}|x'-x''|(|\log 3|+ |\log |x'-x''||)
\\ \nonumber
&&  \qquad\qquad\qquad\quad
\leq 6cc^{(v)}_{\Omega}|\log 3| 2|x'-x''||\log |x'-x''||
\,.
\end{eqnarray} 
Moreover,  
\begin{eqnarray}
\label{prop:v0a4}
\lefteqn{
\int_{\partial\Omega\setminus {\mathbb{B}}_{2}(x',2|x'-x''|)}
 |\, S_{ {\mathbf{a}} }(x'-y)-S_{ {\mathbf{a}} }(x''-y)  \,|
 \,d\sigma_{y}
 }
 \\
\nonumber
&&\qquad\qquad\qquad\qquad\qquad\qquad
\leq
c
\int_{\partial\Omega\setminus {\mathbb{B}}_{2}(x',2|x'-x''|)}
\frac{  |x'-x''|    }{  |x'-y|    }\,d\sigma_{y}
\end{eqnarray} 
Then Lemma 3.5 (iv) of \cite{DoLa17} implies that there exists $ c^{iv}_\Omega\in]0,+\infty[$ such that
 \[
 \int_{\partial\Omega\setminus {\mathbb{B}}_{2}(x',2|x'-x''|)}
\frac{ \,d\sigma_{y}   }{  |x'-y|    } 
\leq    c^{iv}_\Omega|\log|x'-x''|| 
 \]
 for all $x',x''\in \partial\Omega$, $0<|x'-x''|\leq 1/e$.	 
 Hence, the statement holds true.\hfill  $\Box$ 

\vspace{\baselineskip}

Next we prove a regularity statement  for the double layer potential of a constant function. To do so, we need to exploit the tangential derivatives of a function defined on the boundary of an open set of class $C^1$. 
If $l,r\in\{1,\dots,n\}$,  then $M_{lr}$ denotes the tangential derivative 
 operator from $C^{1}(\partial\Omega)$ to $C^{0}(\partial\Omega)$ that takes $f$ to  
 \begin{equation}
\label{mlr}
M_{lr}[f]\equiv \nu_{l}\frac{\partial\tilde{f}}{\partial x_{r}}-
\nu_{r}\frac{\partial\tilde{f}}{\partial x_{l}}\qquad {\text{on}}\ \partial\Omega\,,
\end{equation}
where  $\tilde{f}$ is any continuously differentiable  extension of $f$  to an open neighborhood of $\partial\Omega$. We note that $M_{lr}[f]$ is independent of the specific choice of $\tilde{f}$ (cf.~\textit{e.g.}, 
reference \cite[\S 2.21]{DaLaMu21} with Dalla Riva and Musolino). Then we can state the following.
\begin{lem}\label{lem:dl1sgen}
 Let $n\in {\mathbb{N}}\setminus\{0\}$. Let $\Omega$ be a bounded open subset of ${\mathbb{R}}^{n}$ of class $C^{1,1}$.  Let ${\mathbf{a}}$ be as in (\ref{introd0}), (\ref{ellip}), (\ref{symr}).  Let $S_{ {\mathbf{a}} }$ be a fundamental solution of $P[{\mathbf{a}},D]$. Then $W_\Omega[{\mathbf{a}},S_{{\mathbf{a}}}   ,1]\in C^{1,\omega_1(\cdot)}(\partial\Omega)$.
\end{lem}{\bf Proof.} By reference \cite[Thm.~9.1]{DoLa17} with Dondi, we know that $W_\Omega[{\mathbf{a}},S_{{\mathbf{a}}}   ,1]\in C^{1}(\partial\Omega)$ and that the tangential derivatives of $W_\Omega[{\mathbf{a}},S_{{\mathbf{a}}}   ,1]$ are delivered by the following formula.
 \begin{eqnarray}
\label{lem:dl1sgen1}
\lefteqn{
M_{lj}[W_\Omega[{\mathbf{a}},S_{{\mathbf{a}}}   ,1]] 
}
\\ \nonumber
&& \qquad
=\nu_{l} Q_j\left[\nu\cdot a^{(1)},1\right] 
-\nu_{j} Q_l\left[\nu\cdot a^{(1)},1\right] 
\\ \nonumber
&& \qquad\quad 
+\nu \cdot a^{(1)}
\left\{
Q_l\left[\nu_{j},1\right] 
-
Q_j\left[ \nu_{l},1\right] 
\right\}
+R[\nu_{l},\nu_{j},1]
\qquad{\mathrm{on}}\ \partial\Omega\,,
\end{eqnarray}
where
\[
Q_j[g,\mu](x)
=\int_{\partial\Omega}(g(x)-g(y))\frac{\partial S_{ {\mathbf{a}} }}{\partial x_{j}}(x-y)\mu(y)\,d\sigma_{y}\quad\forall x\in \partial\Omega\,,
\]
for all $(g,\mu)\in  C^{0,1}(\partial\Omega)\times L^{\infty}(\partial\Omega)$ 
    and
\begin{eqnarray*}
\lefteqn{
R[\nu_{l},\nu_{j},1]\equiv \sum_{r=1}a_{r} 
\left\{
Q_r[\nu_{l}\nu_{j},1]-\nu_{l}
Q_r[\nu_{j},1]
-Q_r[\nu_{j},\nu_{l}]
\right\}
}
\\ \nonumber
&&\qquad\qquad\qquad\qquad\qquad\qquad
+a\left\{
\nu_{l}v_\Omega[ S_{ {\mathbf{a}} }   ,\nu_{j}   ]
- \nu_{j} v_\Omega[S_{ {\mathbf{a}} }  ,\nu_{l}   ]
\right\}\qquad\text{on}\ \partial\Omega\,,
\\ \nonumber
\lefteqn{
v_\Omega[ S_{ {\mathbf{a}} }   ,\nu_{j}   ](x)
\equiv\int_{\partial\Omega}S_{ {\mathbf{a}} }(x-y) \nu_{j}(y)\,d\sigma_{y}\qquad\forall x\in {\mathbb{R}}^n 
}
\end{eqnarray*} 
for all $l,j\in\{1,\dots,n\}$.
 By the Lipschitz continuity of the components of $\nu$, Proposition \ref{prop:v0a} above and Theorem 7.2 of \cite{DoLa17}  imply that $v_\Omega[ S_{ {\mathbf{a}} }   ,\nu_{j}   ]$
  belongs to $C^{0,\omega_1(\cdot)}(\partial\Omega)$. By the Lipschitz continuity of the components of $\nu$,  Theorem 8.2 (i) of \cite{DoLa17} implies that
 $Q_r[\nu_{l}\nu_{j},1]$, 
$Q_r[\nu_{j},1]$, $Q_j\left[\nu\cdot a^{(1)},1\right]$, $Q_r[\nu_{j},\nu_{l}]$, belong  to $C^{0,\omega_1(\cdot)}(\partial\Omega)$
 for all $j$, $l$, $r\in\{1,\dots,n\}$. Hence, the tangential derivatives $M_{lj}[W_\Omega[{\mathbf{a}},S_{{\mathbf{a}}}   ,1]] $ belong to $C^{0,\omega_1(\cdot)}(\partial\Omega)$  for all $j$, $l\in\{1,\dots,n\}$, and accordingly $W_\Omega[{\mathbf{a}},S_{{\mathbf{a}}}   ,1]$ belongs to $C^{1,\omega_1(\cdot)}(\partial\Omega)$ (cf.~\textit{e.g.},  \cite[Lem.~2.2]{DoLa17}).  \hfill  $\Box$ 

\vspace{\baselineskip}

 As a consequence of Lemmas \ref{lem:boest}, \ref{lem:tgdlgen}, 
 \ref{lem:dl1sgen}, we can apply Theorem \ref{thm:iokreg} and prove the following theorem on the continuity of the double layer potential on the boundary.
\begin{thm}\label{thm:dllreggen}
  Let  $\beta\in]0,1]$.  Let $\Omega$ be a bounded open subset of ${\mathbb{R}}^{n}$ of class  $C^{1,1}$. 		 
  Assume that  the following condition holds
  \begin{equation}\label{thm:mftgdlgech1}
\sup_{x\in\partial\Omega}\sup_{r\in]0,+\infty[}
\left|
\int_{(\partial\Omega)\setminus  {\mathbb{B}}_n(x,r)}
{\mathrm{grad}}_{\partial\Omega,x}\overline{B^{*}_{\Omega,y}}\left(S_{{\mathbf{a}}}(x-y)\right)
\,d\sigma_y
\right|<+\infty\,,
\end{equation}
\textit{i.e.}, the maximal function of the tangential gradient of the kernel of the double layer potential with respect to its first variable is bounded.	
 
Let ${\mathbf{a}}$ be as in (\ref{introd0}), (\ref{ellip}), (\ref{symr}).  Let $S_{ {\mathbf{a}} }$ be a fundamental solution of $P[{\mathbf{a}},D]$. Then the following statements hold.
  \begin{enumerate}
\item[(i)] If $\beta<1$, then the  operator $W_\Omega[{\mathbf{a}},S_{{\mathbf{a}}}   ,\cdot]$ from 
$C^{0,\beta}(\partial\Omega)$ to $C^{1, \beta }(\partial\Omega)$ defined by 
 (\ref{introd3}) for all $\mu\in C^{0,\beta}(\partial\Omega)$ is linear and continuous.
 \item[(ii)] If $\beta=1$, then the  operator $W_\Omega[{\mathbf{a}},S_{{\mathbf{a}}}   ,\cdot]$ from 
$C^{0,1}(\partial\Omega)$ to $C^{1, \omega_1(\cdot) }(\partial\Omega)$ defined by 
 (\ref{introd3})  for all $\mu\in C^{0,1}(\partial\Omega)$ is linear and continuous.
\end{enumerate}
\end{thm}
{\bf Proof.}   By formula (\ref{eq:boest3}), we have $
\overline{B^{*}_{\Omega,y}}\left(S_{{\mathbf{a}}}(\cdot-y)\right))\in C^1((\partial\Omega)\setminus\{y\})$ for all $y\in \partial\Omega$.	 
If $n=2$, we choose $\epsilon\in]0,1[$ and Lemma  \ref{lem:boest} (ii), (iii) implies that the kernel of the double layer potential belongs to   ${\mathcal{K}}_{\epsilon, 1,1}(\partial\Omega\times\partial\Omega)$. Then the imbedding Proposition \ref{prop:kerem} (ii)  implies that ${\mathcal{K}}_{\epsilon, 1,1}(\partial\Omega\times\partial\Omega)$ is contained in ${\mathcal{K}}_{\epsilon, 1+\epsilon,1}(\partial\Omega\times\partial\Omega)$.

If $n\geq 3$ Lemma  \ref{lem:boest} (i), (iii) implies that the kernel of the double layer potential belongs to the class ${\mathcal{K}}_{n-2, n-1,1}(\partial\Omega\times\partial\Omega)$.

Then if $n\geq 2$ Lemma  \ref{lem:tgdlgen} and 
 condition (\ref{thm:mftgdlgech1})  
  imply that the tangential gradient with respect to the variable $x$ of the kernel of the double layer potential belongs to the class $({\mathcal{K}}_{n-1,n,1}^\sharp(\partial\Omega\times\partial\Omega))^n $. We now plan to apply Theorem \ref{thm:iokreg} (ii). By Lemma \ref{lem:dl1sgen}, we have 
\[
W_\Omega[{\mathbf{a}},S_{{\mathbf{a}}}   ,1]\in C^{1,\omega_1(\cdot)}(\partial\Omega)
\subseteq C^{1,\alpha}(\partial\Omega)\qquad\forall \alpha\in]0,1[\,.
\]
 Moreover,
 \begin{eqnarray*}
&&\beta\leq 1\leq n-1 \equiv t_1 <(n-1)+\beta\,,
\\
&&t_2\equiv n\geq 2>\beta\,,
\quad
s_1\equiv
\left\{
\begin{array}{ll}
\epsilon<2-1 =n-1 &  \text{if}\ n=2\,,
 \\
 (n-1)-1<n-1 & \text{if}\ n\geq 3\,.
\end{array}
\right.
\end{eqnarray*}
(i) If $\beta<1$, then 
\[
t_2-\beta=n-\beta>n-1\,,\quad
t_2=n<(n-1)+\beta+1=(n-1)+\beta+t_3
\ \  \text{where}\ t_3\equiv 1\,.
\]
and $W_\Omega[{\mathbf{a}},S_{{\mathbf{a}}}   ,1]\in C^{1,\omega_1(\cdot)}(\partial\Omega)\subseteq C^{1,\min\{\beta,(n-1)+t_3+\beta-t_2\}}(\partial\Omega)$. Thus Theorem \ref{thm:iokreg} (ii) (b) implies that $W_\Omega[{\mathbf{a}},S_{{\mathbf{a}}}   ,\cdot]$ is linear and continuous from $C^{0,\beta}(\partial\Omega)$ to
\[
C^{1,\min\{\beta,(n-1)+t_3+\beta-t_2\}}(\partial\Omega)
=C^{1,\min\{\beta,(n-1)+1+\beta- n )\}}(\partial\Omega)
=C^{1,\beta}(\partial\Omega)\,.
\]
(ii) If $\beta=1$, then  $t_2-\beta=n-\beta=n-1$ and   
 $W_\Omega[{\mathbf{a}},S_{{\mathbf{a}}}   ,1]\in C^{1,\omega_1(\cdot)}(\partial\Omega)\subseteq C^{1,
\max\{r^{ \beta },\omega_{1}(r)\}
}(\partial\Omega)$. Thus Theorem \ref{thm:iokreg} (ii) (bb) implies that $W_\Omega[{\mathbf{a}},S_{{\mathbf{a}}}   ,\cdot]$ is linear and continuous from $C^{0,\beta}(\partial\Omega)=C^{0,1}(\partial\Omega)$ to
\[
C^{1,
\max\{r^{ \beta },\omega_{1}(r)\}
}(\partial\Omega)=C^{1,
\max\{r^{1 },\omega_{1}(r)\}
}(\partial\Omega)=C^{1,\omega_{1}(\cdot)}(\partial\Omega)\,.
\]
 and thus the proof is complete.\hfill  $\Box$ 

\vspace{\baselineskip}

 For the validity of condition (\ref{thm:mftgdlgech1}), we refer to \cite{La23b}.

\section{Backmatter}


\subsection*{Funding and/or Conflicts of interests/Competing interests}

 The author  acknowledges  the support of the Research 
Project GNAMPA-INdAM   $\text{CUP}\_$E53C22001930001 `Operatori differenziali e integrali in geometria spettrale'  and is indebted to Prof. Otari Chkadua for a number of references. 
Data sharing is not applicable to this article as no data sets were generated or analysed during the current study. This paper does not have any  conflict of interest or competing interest.


\begin{thebibliography}{11}

 \bibitem{Cha84}
 I.~Chavel, {\em  Eigenvalues in Riemannian geometry}.
Including a chapter by Burton Randol. With an appendix by Jozef Dodziuk. Pure and Applied Mathematics, 115. Academic Press, Inc., Orlando, FL, 1984.

 
\bibitem{Chk23}
O.~Chkadua, {\em Personal communication}, 2023.
 

 \bibitem{Ci95}
A.~Cialdea, {\em A general theory of hypersurface potentials}, Annali di Matematica Pura ed Applicata, {\bf 168} (1995), pp 37--61.


\bibitem{CoKr83}
D.~Colton and R.~Kress, {\em Integral equation methods in scattering 
theory},  Wiley, New York, 1983. 

\bibitem{Da13}
M.~Dalla Riva, {\em A family of fundamental solutions of elliptic partial differential operators with real constant coefficients}. Integral Equations Operator Theory, {\bf 76} (2013),  1--23.

\bibitem{DaMoMu13} 
 M.~Dalla Riva, J.~Morais, and P.~Musolino, {\em A family of fundamental solutions of elliptic partial differential operators with quaternion constant coefficients}. Math.~Methods Appl.~Sci.,  {\bf 36} (2013),  1569--1582.


\bibitem{DaLaMu21}
M. Dalla Riva, M.~Lanza de Cristoforis, and P.~Musolino, {\em Singularly Perturbed Boundary Value Problems. A Functional Analytic Approach},  Springer, Cham,  2021.


\bibitem{DoLa17}
F.~Dondi and M.~Lanza de Cristoforis, {\em  Regularizing properties
 of the double layer potential
of  second order elliptic differential operators},   Mem. Differ. Equ.
Math. Phys. 71 (2017), 69--110.

\bibitem{FiDe70}
G.~Fichera, L.~De Vito, {\em Funzioni analitiche di una variabile complessa}. Libreria Eredi Virgilio Veschi, Roma 1970.

 
\bibitem{Ge67}
T.G.~Gegelia, {\em
Certain special classes of functions and their properties.} (Russian), 
 Sakharth. SSR Mecn. Akad. Math. Inst. \v{S}rom. {\bf 32} (1967), 94--139.
 

\bibitem{Gi32}
G.~Giraud, {\em Sur certains problemes non lineaires de Neumann et sur certains problemes non lineaires mixtes}. Ann. Ec. N. Sup. {\bf 49} (1932), 1--104 and 245--308. 

 
\bibitem{Gi34}
G.~Giraud, {\em 
\'{E}quations à int\'{e}grales principales; \'{e}tude suivie d'une application.} (French)
Ann. Sci. \'{E}cole Norm. Sup.  {\bf 51} (1934), 251--372.
 
  
\bibitem{Gu67}
N.M.~G\"{u}nter, {\em
Potential theory and its applications to basic problems of mathematical physics}, 
translated from the Russian by John R. Schulenberger, Frederick Ungar Publishing Co., New York, 1967.

\bibitem{He92}
U.~Heinemann, {\em Die regularisierende Wirkung der Randintegraloperatoren der klassischen Potentialtheorie in den R\"{a}umen h\"{o}lderstetiger Funktionen}, Diplomarbeit, Universit\"{a}t Bayreuth, 1992.

 
\bibitem{HsWe08}
G.C.~Hsiao and W.L.~Wendland,
 {\em Boundary integral equations}, volume 164 of {\em Applied
  Mathematical Sciences}.
 Springer-Verlag, Berlin, 2008. 
 

\bibitem{Jo55}
F.~John,  {\em Plane waves and spherical means applied to partial differential
  equations}.  Interscience Publishers, New York-London, 1955.
 
 \bibitem{Ki89}
A.~Kirsch, {\em Surface gradients and continuity properties for some 
integral operators in classical scattering theory}, Math.~Methods  
Appl.~Sciences, {\bf 11} (1989), 789--804. 


\bibitem{KiHe15}
A.~Kirsch and  F.~Hettlich, {\em  The Mathematical Theory of Time-Harmonic Maxwell's Equations; Expansion-, Integral-, and Variational Methods},  Springer, 2015.  

\bibitem{KuGeBaBu79}
V.D.~Kupradze, T.G.~Gegelia, M.O.~Basheleishvili and 
T.V.~Burchuladze, {\em Three-dimensional Problems of the Mathematical 
Theory of Elasticity and Thermoelasticity}, North-Holland Publ.~Co., Amsterdam, 
1979. 

\bibitem{La20}  Lanza de Cristoforis M.~{\em An inequality for {H}\"{o}lder continuous functions generalizing a
  result of Carlo Miranda},  Computational and Analytic Methods in Science and Engineering,
C. Constanda (ed.), Birkh\"auser-Springer, New York, 2020, pp. 197--221. 

\bibitem{La22b}
M.~Lanza de Cristoforis, {\em Classes of kernels and continuity properties of the tangential gradient of an integral operator in H\"{o}lder spaces on a manifold}, submitted, (2022).

 
\bibitem{La23b}
 M.~Lanza de Cristoforis, {\em On the tangential gradient of the kernel of the double layer potential}, submitted, (2023).
 
 

 \bibitem{Mik70}
S.G.~Mikhlin, {\em Mathematical physics, an advanced course},  translated from the Russian, North-Holland Publishing Co., Amsterdam-London,  1970.

 
 \bibitem{MikPr86}
S.G.~Mikhlin and S.~Pr\"{o}ssdorf, {\em Singular integral  Operators}, Springer-Verlag, Belin, 1986. 
 
 
\bibitem{Mi65}
C.~Miranda, {\em Sulle propriet\`{a} di regolarit\`{a} di certe 
trasformazioni integrali}, Atti Accad. Naz.   
Lincei Mem. Cl. Sci. Fis. Mat. Natur. Sez I, {\bf 7} (1965), 303--336. 

\bibitem{Mi70} 
C.~Miranda, {\em Partial differential equations of elliptic type, 
Second revised edition},   Springer-Verlag, Berlin,  1970.

\bibitem{Mit14}
M.~Mitrea, {\em The almighty double layer: recent perspectives}. Invited presentation at the 13$th$ International Conference on Integral methods in Science and Engineering, July 21--25, 2014. 

 
\bibitem{MitMit13}
I.~Mitrea and M.~Mitrea, {\em Multi-Layer Potentials and Boundary Problems, for Higher-Order Elliptic Systems in Lipschitz Domains}.
Lecture Notes in Mathematics, Springer, Berlin, {\it etc.} 2013.

 

\bibitem{Sc31}
J.~Schauder,  {\em Potentialtheoretische Untersuchungen}, Math. 
Z., {\bf 33}  (1931),  602--640.

\bibitem{Sc32}
J.~Schauder,  {\em Bemerkung zu meiner Arbeit ``Potentialtheoretische Untersuchungen I (Anhang)''}, Math.  Z., {\bf 35} (1932), 536--538.

 \end{thebibliography}
\end{document}